\documentclass[12pt,reqno]{amsart}
\usepackage{amsmath}
\usepackage{amssymb}
\usepackage{verbatim}
\usepackage{mathrsfs}
\usepackage[capitalise]{cleveref}
\usepackage[noadjust]{cite}

\usepackage[top=0.68in,bottom=0.68in,left=0.7in,right=0.7in]{geometry}

\newcommand{\C}{\mathbb{C}}    
\newcommand{\N}{\mathbb{N}}    
\newcommand{\NN}{\mathbb{N}_0} 
\newcommand{\R}{\mathbb{R}}    
\newcommand{\Z}{\mathbb{Z}}    
\newcommand{\sD}{\mathsf{D}}
\newcommand{\CH}[1]{\mathscr{C}^{#1}(\R)} \newcommand{\Lp}[1]{L_{#1}(\mathbb{R})}
\newcommand{\cd}{\mathcal{R}}  
\newcommand{\sm}{\operatorname{sm}}  

\newcommand{\bp}{ \begin{proof} }
\newcommand{\ep}{\hfill \end{proof} }
\newcommand{\be}{ \begin{equation} }
\newcommand{\ee}{ \end{equation} }
\newcommand{\imply}{ \Longrightarrow }
\newcommand{\tp}{\mathsf{T}}  
\newcommand{\fs}{\operatorname{fsupp}}
\newcommand{\supp}{\operatorname{supp}} 
\newcommand{\mspan}{\operatorname{span}}

\newcommand{\pp}{\mathsf{p}}
\newcommand{\pq}{\mathsf{q}}
\newcommand{\lp}[1]{l_{#1}(\mathbb{Z})}
\newcommand{\lrs}[3]{(l_{#1}(\mathbb{Z}))^{#2\times #3}}

\newcommand{\PL}{\Pi}   
\newcommand{\PR}{\mathscr{P}}   
\newcommand{\PV}{\mathscr{V}}
\newcommand{\PB}{\mathscr{B}}

\newcommand{\sd}{\mathcal{S}}  

\newcommand{\wh}{\widehat}
\renewcommand{\le}{\leqslant}
\renewcommand{\ge}{\geqslant}
\newcommand{\bs}{\backslash}
\newcommand{\ol}{\overline}

\newcommand{\bo}{\mathscr{O}} 
\newcommand{\setsp}{\;:\;}     

\newcommand{\vgu}{\upsilon} 
\newcommand{\td}{\boldsymbol{\delta}}

\newtheorem{lemma}{Lemma}
\newtheorem{prop}[lemma]{Proposition}
\newtheorem{cor}[lemma]{Corollary}
\newtheorem{theorem}[lemma]{Theorem}

\newtheorem{example}{Example}

\newtheorem{definition}[lemma]{Definition}

\numberwithin{equation}{section}
\numberwithin{lemma}{section}

\begin{document}

\title{Analysis and Convergence of Hermite Subdivision Schemes}

\author{Bin Han}

\address{Department of Mathematical and Statistical Sciences,
University of Alberta, Edmonton,\quad Alberta, Canada T6G 2G1.
\quad {\tt bhan@ualberta.ca}\quad {\tt http://www.ualberta.ca/$\sim$bhan}
}

\thanks{Research was supported in part by the Natural Sciences and Engineering Research Council of Canada (NSERC).
The author thanks Professor Nira Dyn for discussing Hermite subdivision schemes.
}


\makeatletter \@addtoreset{equation}{section} \makeatother

\begin{abstract}
Hermite interpolation property is desired in applied and computational mathematics. Hermite and vector subdivision schemes are of interest in CAGD for generating subdivision curves and in computational mathematics for building Hermite wavelets to numerically solve partial differential equations. In contrast to well-studied scalar subdivision schemes, Hermite and vector subdivision schemes employ matrix-valued masks and vector input data, which make their analysis much more complicated and difficult than their scalar counterparts. Despite recent progresses on Hermite subdivision schemes, several key questions still remain unsolved, for example, characterization of Hermite masks, factorization of matrix-valued masks, and convergence of Hermite subdivision schemes. In this paper, we shall study
Hermite subdivision schemes through investigating vector subdivision operators acting on vector polynomials and establishing the relations among Hermite subdivision schemes, vector cascade algorithms and refinable vector functions. This approach allows us to resolve several key problems on Hermite subdivision schemes including characterization of Hermite masks, factorization of matrix-valued masks, and convergence of Hermite subdivision schemes.
\end{abstract}


\maketitle

\pagenumbering{arabic}

\section{Introduction, Motivations and Main Results}

A vector subdivision scheme is an iterative averaging algorithm by recursively applying a vector subdivision operator to a given initial input vector sequence.
A Hermite subdivision scheme is a special modified type of vector subdivision schemes for computing a limiting function and its consecutive derivatives.
Due to their highly desired properties such as interpolation, smoothness and short support of basis functions,
Hermite subdivision schemes are of particular interest and importance in CAGD for generating subdivision curves and in computational mathematics for constructing Hermite wavelets to numerically solve partial differential equations.

Let us first recall some notations and definitions related to a Hermite subdivision scheme. Let $r,s\in \N$ be positive integers.
By $\lrs{}{s}{r}$ we denote the linear space of all sequences $u=\{u(k)\}_{k\in \Z}: \Z \rightarrow \C^{s\times r}$.
Similarly,
$\lrs{0}{s}{r}$ consists of all finitely supported sequences $u\in \lrs{}{s}{r}$ with $\{k\in \Z \setsp u(k)\ne 0\}$ being finite.
Let $a=\{a(k)\}_{k\in \Z}\in \lrs{0}{r}{r}$ with $a(k)\in \C^{r\times r}$ for all $k\in \Z$, which is often called a matrix-valued \emph{mask} in CAGD and a matrix-valued \emph{filter} in wavelet theory. A vector or Hermite subdivision scheme can be conveniently expressed through the vector subdivision operator $\sd_a: \lrs{}{s}{r} \rightarrow \lrs{}{s}{r}$ which is defined to be
\be \label{sd}
(\sd_a v)(j):=2\sum_{k\in \Z} v(k) a(j-2k),\qquad j\in \Z
\ee
for $v=\{v(k)\}_{k\in \Z}\in \lrs{}{s}{r}$. In computational mathematics and CAGD, one is often interested in the special case $s=1$ and real-valued sequences/masks $a$. Because complex-valued masks are of interest in certain applications (e.g., see \cite{hanbook}), without restricting ourselves, we deal with masks of
both complex and real numbers. For $u\in \lrs{0}{s}{r}$, we define its \emph{symbol} or \emph{Fourier series} to be
\be \label{fourier}
\wh{u}(\xi):=\sum_{k\in \Z} u(k) e^{-ik\xi},\qquad \xi\in \R,
\ee
which is an $s\times r$ matrix of $2\pi$-periodic trigonometric polynomials.
For $v\in \lrs{0}{s}{r}$, it follows directly from the definition of a vector subdivision operator in \eqref{sd} that
$\wh{\sd_a v}(\xi)=2\wh{v}(2\xi)\wh{a}(\xi)$.

We now recall a Hermite subdivision scheme. Let $w_0: \Z \rightarrow \C^{1\times r}$ be a sequence on $\Z$ (which is a row vector sequence) standing for a given input vector sequence/data.
In a Hermite subdivision scheme of order $r$ (i.e., degree $r-1$) associated with a mask $a\in \lrs{0}{r}{r}$, a sequence of Hermite refinements $w_n: \Z \rightarrow \C^{1\times r}$ for $n\in \N$ is obtained through recursively applying the vector subdivision operator $\sd_a$ on $w_0$ iteratively as follows:
\be \label{hsd:wn}
w_{n}:=(\sd_a^n w_0) \sD^{-n},\qquad n\in \N \quad \mbox{with}\quad
\sD:=\mbox{diag}(1, 2^{-1},\ldots, 2^{1-r}).
\ee
In terms of Fourier series for $w_0\in \lrs{0}{1}{r}$, we have
\be \label{hsd:wn:2}
\wh{w_n}(\xi)=2^n \wh{w_0}(2^n\xi)\wh{a_n}(\xi) \sD^{-n}\quad \mbox{with}\quad
\wh{a_n}(\xi):=\wh{a}(2^{n-1}\xi)\cdots \wh{a}(2\xi)\wh{a}(\xi).
\ee
If $w_{n+1}(2k)=w_n(k)$ for all $k\in \Z$ and $n\in \N$, then such a Hermite subdivision scheme of order $r$ is often called \emph{an interpolatory Hermite subdivision scheme of order $r$}. It is easy to deduce from \eqref{sd} that $w_{n+1}(2k)=w_n(k)$ for all $k\in \Z$ and $n\in \N$ if and only if
\be \label{sd:int}
a(0)=\mbox{diag}(2^{-1},2^{-2},\ldots,2^{-r}) \quad \mbox{and}\quad a(2k)=0,\qquad \forall\; k\in \Z\bs\{0\}.
\ee
That is, an interpolatory Hermite subdivision scheme of order $r$ associated with mask $a\in \lrs{0}{r}{r}$ is just a Hermite subdivision scheme whose mask satisfies the interpolation condition in \eqref{sd:int}. Its underlying basis vector function $h=(h_1,\ldots,h_r)^\tp$ associated with an interpolatory Hermite subdivision scheme of order $r$ often is \emph{a Hermite interpolant of order $r$}, that is, $h\in (\CH{r-1})^r$ and
\be \label{hermite:int}
[h,h',\ldots,h^{(r-1)}](k)=\td(k)I_r,\qquad \forall\, k\in \Z,
\ee
where $\td$ is the \emph{Dirac sequence} such that $\td(0):=1$ and $\td(k)=0$ for all $k\in \Z \bs \{0\}$. Note that $\wh{\td}(\xi)=1$.

The study of interpolatory Hermite subdivision schemes of order $2$ (i.e., degree $1$) has been initiated by Merrien \cite{mer92} and Dyn and Levin \cite{dl95}.
For a vector subdivision scheme,
the matrix $\sD^{-n}$ in \eqref{hsd:wn} and \eqref{hsd:wn:2} is removed for refinements. Vector and Hermite subdivision schemes have been extensively studied in the literature by a lot of researchers, e.g., see \cite{ch19,ccs16,cmss19,dm09,dl95,dl99,han01,han03,hanbook,hj98,hyx05,hz09,jj02,jrz98,mer92,ms11,ms17,ms19,ms98,zhou00} and many references therein
from different perspectives and purposes. We mention that multivariate refinable Hermite interpolants have been characterized in \cite[Corollary~5.2]{han03} (also see \cite[Theorem~6.2.3]{hanbook}) and multivariate noninterpolatory Hermite subdivision schemes have been studied in \cite{hyx05}.

To avoid potential confusion about notations used in this paper, it is important to pay attention to the notation differences of a Hermite subdivision scheme defined here in comparison with other papers in the literature.
The vector subdivision operator in \eqref{sd} and Hermite subdivision schemes in \eqref{hsd:wn} follow the same notations as in the book \cite{hanbook} on framelets and wavelets. As we shall see later in this paper, these notations allow us conveniently link Hermite subdivision schemes to vector cascade algorithms and refinable vector functions in wavelet theory. More precisely,
the masks $\mathbf{A}$ and generated Hermite refinement data $\{f_n\}_{n=1}^\infty$ (or $\{\mathbf{c}_n\}_{n=1}^\infty$) in other papers such as \cite{ch19,ccs16,dl95,dl99,mer92,ms11,ms17,ms19} correspond to $2a^\tp$ for masks and $\{w_n^\tp\}_{n=1}^\infty$ for Hermite refinement data in this paper. In short, the main difference of notations is a transpose put on masks and refinement data.

For $\ell=1,\ldots,r$,
by $e_{\ell}\in \R^r$ we denote the standard unit coordinate (column) vector such that $e_{\ell}=(0,\ldots,0,1,0,\ldots,0)^\tp$ with the only nonzero coefficient $1$ at the $\ell$th entry.
To explain our motivations of this paper, we recall the definition of a convergent Hermite subdivision scheme, e.g., see \cite[Definition~1]{dm09},
\cite[Definition~1.1]{hyx05}, \cite[Definition~2]{ms19} and references therein.

\begin{definition}\label{def:hsd}
{\rm
Let $r\in \N$ and $m\in \NN:=\N\cup\{0\}$ with $m\ge r-1$.
We say that \emph{a Hermite subdivision scheme of order $r$ (i.e., degree $r-1$) associated with a finitely supported matrix-valued mask $a\in \lrs{0}{r}{r}$
is convergent with limiting functions in $\CH{m}$} if for every input vector sequence $w_0=\{w_0(k)\}_{k\in \Z}: \Z \rightarrow \C^{1\times r}$ (i.e., $w_0\in \lrs{}{1}{r}$), there exists a $\CH{m}$ function $\eta: \R\rightarrow \C$ such that for all constants $K>0$,
\be \label{hsd:converg}
\lim_{n\to \infty} \max_{k\in \Z\cap [-2^nK, 2^n K]} \|w_n(k)-[\eta(2^{-n}k), \eta'(2^{-n}k),\ldots, \eta^{(r-1)}(2^{-n}k)]\|_{\ell_\infty}=0.
\ee
Or equivalently,
\be \label{hsd:converg:2}
\lim_{n\to \infty} \max_{k\in \Z\cap [-2^nK, 2^n K]} |w_n(k) e_{\ell+1}- \eta^{(\ell)}(2^{-n}k)|=0,\qquad \forall\; \ell=0,\ldots,r-1,
\ee
where $\eta^{(\ell)}$ stands for the $\ell$th classical derivative of the function $\eta\in \CH{m}$.
}
\end{definition}

Let $r\in \N$ and $m\in \NN$ with $m\ge r-1$. In the following, we first discuss the underlying basis vector function $\phi$ in a convergent
Hermite subdivision scheme of order $r$ associated with a matrix-valued mask $a\in \lrs{0}{r}{r}$ with limiting functions in $\CH{m}$.
For each $\ell=1,\ldots,r$, take $w_0=\td e_\ell^\tp\in \lrs{0}{1}{r}$ (i.e., $w_0(0)=e_\ell^\tp$ and $w_0(k)=(0,\ldots,0)$ for all $k\in \Z\bs\{0\}$) as an initial input vector sequence. Then recursively compute Hermite refinements $\{w_n\}_{n=1}^\infty$ as in \eqref{hsd:wn}.
We define the filter support of the mask $a\in \lrs{0}{r}{r}$  to be $\fs(a):=[L,R]$ with $L,R\in \Z$ such that $a(L)\ne 0$, $a(R)\ne 0$ and $a(k)=0$ for all $k\in \Z\bs [L,R]$.
For simplicity of discussion, we assume that $0\in \fs(a)$. Using induction on $n$, one can easily observe that $\fs(w_n)\subseteq (2^n-1)\fs(a)$ for all $n\in \N$.
By \cref{def:hsd}, since the mask $a$ has finite support, there exists a compactly supported $\CH{m}$ function $\phi_\ell: \R \rightarrow \C$ such that \eqref{hsd:converg:2} holds with $\eta=\phi_\ell$ and $\supp(\phi_\ell)\subseteq \fs(a)$. More explicitly, the basis vector function $\phi:=(\phi_1,\ldots,\phi_r)^\tp$ of the convergent Hermite subdivision scheme is given as the limiting vector function from the initial sequence $w_0=\td I_r$ satisfying
\be \label{sdn:phi}
\lim_{n\to \infty}
\|(\sd_a^n (\td I_r))(\cdot)\sD^{-n}-[\phi,\phi',\ldots,\phi^{(r-1)}](2^{-n}\cdot) \|_{\lrs{\infty}{r}{r}}=0
\quad \mbox{with}\quad
\phi:=(\phi_1,\ldots,\phi_{r})^\tp,
\ee
where $I_r$ stands for the $r\times r$ identity matrix and $\td$ is the Dirac sequence. Note that
\be \label{an}
\sd_a^n (\td I_r)=2^n a_n
\quad \mbox{with}\quad
\wh{a_n}(\xi):=\wh{a}(2^{n-1}\xi)\cdots \wh{a}(2\xi)\wh{a}(\xi).
\ee
In terms of entries,
\eqref{sdn:phi} can be equivalently rewritten as
\be \label{sdn:phi:2}
\lim_{n\to \infty}
\|2^{(\ell-1)n} a_n(\cdot)e_{\ell+1}-\phi^{(\ell)}(2^{-n}\cdot) \|_{\lrs{\infty}{1}{r}}=0,\qquad \forall\, \ell=0,\ldots,r-1.
\ee
Conversely, if there exists a compactly supported vector function $\phi\in (\CH{m})^r$ such that \eqref{sdn:phi} holds, then
the Hermite subdivision scheme of order $r$ associated with mask $a\in \lrs{0}{r}{r}$ must be convergent with limiting functions in $\CH{m}$. Indeed,
for any initial vector sequence $w_0\in \lrs{}{1}{r}$, define refinement data $\{w_n\}_{n=1}^\infty$ as in \eqref{hsd:wn}. Define $\eta:=w_0*\phi:=\sum_{k\in \Z} w_0(k)\phi(\cdot-k)$.
Then we can directly deduce from \eqref{hsd:wn:2} that \eqref{sdn:phi} implies \eqref{hsd:converg}.
Hence, a Hermite subdivision scheme of order $r$ associated with mask $a\in \lrs{0}{r}{r}$ is convergent with limiting functions in $\CH{m}$ if and only if \eqref{sdn:phi} (or equivalently \eqref{sdn:phi:2}) holds for a compactly supported vector function $\phi\in (\CH{m})^r$.

Note that the matrix $\sD^{-n}$ in \eqref{hsd:wn} and \eqref{hsd:wn:2} is removed for refinements in a vector subdivision scheme.
Because vector and Hermite subdivision schemes employ matrix-valued masks and vector sequence data, it is widely known that their analysis is often much more difficult and complicated than their scalar counterparts (e.g., see \cite{cdm91,dl02,hanbook}).
Due to the extra factor $\sD^{-n}$ in \eqref{hsd:wn} and \eqref{hsd:wn:2} for Hermite subdivision schemes, it is known in the literature (e.g., \cite{dm09,ms19}) that a Hermite subdivision scheme is not stationary, i.e., the Hermite refinements $w_n$ in \eqref{hsd:wn} can be rewritten as $w_n=\sd_{\sD^{n-1} a \sD^{-n}} w_{n-1}$, whose masks $\sD^{n-1} a\sD^{-n}$ depend on the subdivision level $n$.
This non-stationarity of level dependent masks makes it much harder to analyze Hermite subdivision schemes than vector subdivision schemes. To study the convergence of a Hermite or vector subdivision scheme, the dominating approach in the area of subdivision schemes is to factorize a Hermite or vector subdivision scheme into another derived associated vector subdivision scheme through the factorization of its underlying masks, e.g., see \cite{ccs16,dm09,dl95,dl99,ms11,ms17,ms19} and references therein for more details.
This approach requires factorization of the symbol of a matrix-valued mask. This is a highly nontrivial task, because factorizing a matrix of $2\pi$-periodic trigonometric polynomials for vector or Hermite subdivision schemes is much harder and complicated than factorizing a $2\pi$-periodic trigonometric polynomial for scalar subdivision schemes.
As a consequence, in recent years, we see a growing interest in the literature to study factorization of the symbol of a matrix-valued mask for a Hermite subdivision scheme, to only mention a few references here, see \cite{ccs16,dm09,ms11,ms17,ms19} and references therein.
Despite recent progresses on Hermite subdivision schemes, a few key questions on Hermite subdivision schemes still remain open, for example,
characterization and factorization of matrix-valued masks, and convergence of a Hermite subdivision scheme.
It is the purpose of this paper to study and resolve these key questions on Hermite subdivision schemes by linking these problems with vector cascade algorithms and refinable vector functions in wavelet theory.

To study convergent Hermite subdivision schemes, it is crucial to characterize their masks with desired properties for convergent Hermite subdivision schemes. This topic is closely
related to how a vector subdivision operator acts on vector polynomials, which in turn is linked to the notion of sum rules of a matrix-valued mask (e.g., see \cite{dl02,han01,han03,hanbook,jp01,jj02}). For a matrix-valued mask $a\in \lrs{0}{r}{r}$ and $m\in \NN:=\N\cup\{0\}$, we way that $a$ has \emph{order $m+1$ sum rules with respect to a (moment) matching filter $\vgu_a\in \lrs{0}{1}{r}$} if $\wh{\vgu_a}(0)\ne 0$ and
\be \label{sr}
\wh{\vgu_a}(2\xi)\wh{a}(\xi)= \wh{\vgu_a}(\xi)+\bo(|\xi|^{m+1})\quad \mbox{and}\quad
\wh{\vgu_a}(2\xi)\wh{a}(\xi+\pi)=\bo(|\xi|^{m+1}),\quad \xi\to 0.
\ee
Here the notation $\wh{f}(\xi)=\wh{g}(\xi)+\bo(|\xi|^{m+1})$ as $\xi \to 0$ simply stands for $\wh{f}^{(j)}(0)=\wh{g}^{(j)}(0)$ for all $j=0,\ldots,m$, as long as both $\wh{f}$ and $\wh{g}$ are smooth near the origin.
In this paper we shall frequently use this short-hand notation, because it makes the presentation much compact and simple.

To study Hermite subdivision schemes,
we shall extensively study in \cref{sec:sd} how a vector subdivision operator acts on vector polynomials.
For $m\in \NN$, by $\PL_m$ we denote the space of all polynomials of degree no more than $m$. An element in $(\PL_m)^{1\times r}$ is simply a row vector of polynomials (i.e., a vector polynomial) having degree no more than $m$. For $\vec{\pp}\in (\PL_m)^{1\times r}$, note that the vector polynomial $\vec{\pp}$ is uniquely determined by the vector polynomial sequence $\{\vec{\pp}(k)\}_{k\in \Z}$ by restricting the vector polynomial $\vec{\pp}$ on $\Z$. Therefore, throughout this paper we do not distinguish a vector polynomial and its corresponding vector polynomial sequence.
We caution the reader that we cannot simply replace the polynomial space $\PL_m$ for scalar subdivision schemes by $(\PL_m)^{1\times r}$ for vector or Hermite subdivision schemes, largely because $\dim(\PL_m)=m+1<r(m+1)=\dim((\PL_m)^{1\times r})$ for $r>1$.
For $m\in \NN$ and $v\in \lrs{0}{1}{r}$, as in \cite{han01,han03,hanbook}, we define the following vector polynomial subspace of $(\PL_m)^{1\times r}$ by
\be \label{Pvm}
\PR_{m,v}:=\{\pp*v \setsp \pp\in \PL_m\} \quad \mbox{with}\quad
\pp*v:=\sum_{k\in \Z} \pp(\cdot-k)v(k).
\ee
Note that $\PR_{m,v}\subseteq (\PL_m)^{1\times r}$.
If $\wh{v}(0)\ne 0$, then $\dim(\PR_{m,v})=\dim(\PL_m)$. The vector polynomial space $\PR_{m,v}$ plays the role for vector and Hermite subdivision schemes as $\PL_m$ for scalar subdivision schemes. For scalar subdivision schemes with $r=1$, we always have $\PR_{m,v}=\PL_m$ for any choice of $v\in \lp{0}$ with $\wh{v}(0)\ne 0$. However, for vector and Hermite subdivision schemes with $r>1$, different choices of $v\in \lrs{0}{1}{r}$ may lead to different subspaces of $(\PL_m)^{1\times r}$ and how to find suitable $v\in \lrs{0}{1}{r}$ for a given matrix-valued mask becomes nontrivial at all.

For a vector subdivision scheme acting on vector polynomials, we shall prove in \cref{sec:sd} the following result, which plays a key role in our study of Hermite subdivision schemes.

\begin{theorem}\label{thm:sd:sr}
Let $a\in \lrs{0}{r}{r}$ and $v\in \lrs{0}{1}{r}$ with $\wh{v}(0)\ne 0$.
For a nonnegative integer $m\in \NN$,
the following statements are equivalent to each other:
\begin{enumerate}
\item[(1)] $\sd_a \PR_{m,v}=\PR_{m,v}$ and $\sd_a (\wh{v}(0))=\wh{v}(0)$, where $\wh{v}(0)$ is regarded as a constant vector sequence.
\item[(2)] $\sd_a \vec{\pp}\in \PR_{\vec{\pp}}:=\mspan\{\vec{\pp},\vec{\pp}',\ldots,\vec{\pp}^{(m)}\}$
    and $\sd_a (\vec{\pp}^{(m)})=\vec{\pp}^{(m)}$, where $\vec{\pp}:=\frac{(\cdot)^m}{m!}*v$.
\item[(3)] There exists a finitely supported sequence $c\in \lp{0}$ with $\wh{c}(0)=1$ such that
\be \label{sr:1}
\wh{v}(2\xi)\wh{a}(\xi)=\wh{c}(\xi) \wh{v}(\xi)+\bo(|\xi|^{m+1})\quad \mbox{and}\quad
\wh{v}(2\xi)\wh{a}(\xi+\pi)=\bo(|\xi|^{m+1}),\quad \xi\to 0.
\ee
\item[(4)] There exists a finitely supported sequence $d\in \lp{0}$ with $\wh{d}(0)=1$ such that the sequence $\vgu_a\in \lrs{0}{r}{r}$ defined by $\wh{\vgu_a}(\xi):=\wh{d}(\xi)\wh{v}(\xi)$ satisfies \eqref{sr} for order $m+1$ sum rules.
%
%
\item[(5)]
$\sd_a\PR_{m,v}\subseteq \PR_{m,v}$ and all the eigenvalues of $\sd_a: \PR_{m,v}\rightarrow \PR_{m,v}$ are $2^{-j}$ for $j=0,\ldots,m$. More precisely,
    $\sd_a \vec{\pp}_j=2^{-j} \vec{\pp}_j$ and $\vec{\pp}_j=\vec{\pp}_m^{(m-j)}$ for $j=0,\ldots,m$, where

\be \label{sd:eigpoly}
\vec{\pp}_j(x):=\left(\frac{(\cdot)^{j}}{j!}*\vgu_a\right)(x)=
\sum_{k=0}^{j} \frac{(-i)^k}{k!(j-k)!} x^{j-k} \wh{\vgu_a}^{(k)}(0)
\ee
and $\vgu_a$ is given in item (4). Note that $\PR_{m,\vgu_a}=\PR_{m,v}$ and $\deg(\vec{\pp}_j)=j$ for all $j=0,\ldots,m$.
\item[(6)] $\sd_a \vec{\pp}_m=2^{-m} \vec{\pp}_m$ for some $\vec{\pp}_m\in \PR_{m,v}$ satisfying $\deg(\vec{\pp}_m)=m$.
\end{enumerate}
Moreover, any of the above items (1)--(6) implies $\sd_a(\pq*\vgu_a)=(\pq(2^{-1}\cdot))*\vgu_a$ for all $\pq\in \PL_m$ and
\be \label{sd:deriv}
\sd_a (\vec{\pq}^{(j)})=2^j [\sd_a \vec{\pq}]^{(j)}
=(\vec{\pq}^{(j)}(2^{-1}\cdot))*a
=\sum_{k=0}^\infty \frac{(-i)^k}{k!2^k} \vec{\pq}^{(j+k)}(2^{-1}\cdot)\wh{a}^{(k)}(0),
\quad \forall\; \vec{\pq}\in \PR_{m,v}, j\in \NN.
\ee
\end{theorem}

For $f\in \Lp{1}$, we define its Fourier transform to be $\wh{f}(\xi):=\int_\R f(x) e^{-i\xi x} dx$ for $\xi\in \R$. The Fourier transform can be naturally extended to tempered distributions through duality.
Now using \cref{thm:sd:sr} and results in \cite{han03,hanbook} on vector cascade algorithms and refinable vector functions, we can characterize masks for all Hermite subdivision schemes as follows.

\begin{theorem}\label{thm:hsd:mask}
Let $r\in \N$ and $m\in \NN$ with $m\ge r-1$.
Assume that the Hermite subdivision scheme of order $r$ associated with a matrix-valued mask $a\in \lrs{0}{r}{r}$ is convergent with limiting functions in $\CH{m}$.
Let the compactly supported vector function $\phi=(\phi_1,\ldots,\phi_r)^\tp \in (\CH{m})^r$ be defined as the limiting basis vector function through \eqref{sdn:phi} using the initial data $w_0=\td I_r\in \lrs{0}{r}{r}$.
If
\be \label{cond:phi}
\mspan\{\wh{\phi}(2\pi k) \setsp k\in \Z\}=\C^r \quad \mbox{and}\quad
\mspan\{\wh{\phi}(\pi+2\pi k) \setsp k\in \Z\}=\C^r,
\ee
then the following statements hold:
\begin{enumerate}
\item[(1)]
$1$ is a simple eigenvalue of $\wh{a}(0)$ and all other eigenvalues of $\wh{a}(0)$ are less than $2^{-m}$ in modulus.
\item[(2)] The mask $a$ must be a Hermite mask of accuracy order $m+1$, that is,
the mask $a$ has order $m+1$ sum rules with respect to some $\vgu_a\in \lrs{0}{1}{r}$ in \eqref{sr} such that
\be \label{vgu:hermite:0}
\wh{\vgu_a}(\xi)=\left(1+\bo(|\xi|),i\xi +\bo(|\xi|^2),\ldots,(i\xi)^{r-1}+
\bo(|\xi|^r)\right),
\qquad \xi\to 0,
\ee
that is, for some $c_1,\ldots,c_r\in \lp{0}$ satisfying $\wh{c_1}(0)=\wh{c_2}(0)=\cdots=\wh{c_r}(0)=1$,
\be \label{vgu:hermite}
\wh{\vgu_a}(\xi)=\left(\wh{c_1}(\xi),i\xi \wh{c_2}(\xi),\ldots,(i\xi)^{r-1}\wh{c_r}(\xi)\right)+\bo(|\xi|^{m+1}),
\qquad \xi\to 0
\ee

\item[(3)] $\phi=2\sum_{k\in \Z} a(k) \phi(2\cdot-k)$ with $e_1^\tp \wh{\phi}(0)=1$ and $\wh{\phi}(\xi)=\lim_{n\to \infty} \big[\prod_{j=1}^n \wh{a}(2^{-j}\xi)\big ]e_1$ for $\xi\in \R$.
\end{enumerate}
\end{theorem}

We shall prove \cref{thm:hsd:mask} in \cref{sec:hsd}.
As we shall explain in \cref{sec:hsd}, the condition in \eqref{cond:phi} allows us to avoid degenerate Hermite subdivision schemes and is quite weak in general.
Also, we shall see in \cref{sec:hsd} that \cref{thm:hsd:mask} allows us to easily construct all Hermite masks of accuracy order $m+1$ through solving a system of equations.

To link Hermite subdivision schemes with vector cascade algorithms, let us recall the vector cascade algorithm. For $1\le p\le \infty$, the \emph{refinement operator}
$\cd_a: (\Lp{p})^r \rightarrow (\Lp{p})^r$
is defined to be
\be \label{cd}
\cd_a f:=2\sum_{k\in \Z} a(k) f(2\cdot-k),\qquad f\in (\Lp{p})^r.
\ee
Built on \cref{thm:hsd:mask} and results in \cite{han03,hanbook},
the following result characterizes convergence of Hermite subdivision schemes and will be proved in \cref{sec:converg}.

\begin{theorem}\label{thm:hsd:converg}
Let $r\in \N$ and $m\in \NN$ with $m\ge r-1$.
Let $a\in \lrs{0}{r}{r}$ be a finitely supported mask such that
items (1) and (2) of \cref{thm:hsd:mask}
are satisfied.
Let $\phi=(\phi_1,\ldots,\phi_r)^\tp$ be the vector of compactly supported distributions satisfying $\wh{\phi}(2\xi)=\wh{a}(\xi)\wh{\phi}(\xi)$ and $\wh{\vgu_a}(0) \wh{\phi}(0)=1$.
If the integer shifts of $\phi$ are stable, i.e.,
\be \label{stability}
\mspan\{\wh{\phi}(\xi+2\pi k) \setsp k\in \Z\}=\C^r, \qquad \forall\, \xi\in \R,
\ee
then the following statements are equivalent to each other:
\begin{enumerate}
\item[(1)] The Hermite subdivision scheme associated with the mask $a\in \lrs{0}{r}{r}$ is convergent with limiting functions in $\CH{m}$.
\item[(2)] The vector cascade algorithm associated with the mask $a$ is convergent in $\CH{m}$, that is, for every compactly supported vector function $f\in (\CH{m})^r$ such that
\be \label{initialf}
\wh{\vgu_a}(0)\wh{f}(0)=1 \quad \mbox{and}\quad
\wh{\vgu_a}(\xi)\wh{f}(\xi+2\pi k)=\bo(|\xi|^{m+1}),\quad \xi \to 0, \forall\; k\in \Z\bs\{0\},
\ee
the cascade sequence $\{\cd_a^n f\}_{n=1}^\infty$ converges to $\phi$ in $\CH{m}$, i.e.,
\be \label{cd:converg}
\phi \in \CH{m} \quad \mbox{and}\quad \lim_{n\to \infty} \|\cd_a^n f-\phi\|_{(\CH{m})^r}=0.
\ee
\item[(3)] $\sm_\infty(a)>m$, where the quantity $\sm_\infty(a)$ is stated in \cref{sec:converg} and is introduced in \cite{han03,hanbook}.
\end{enumerate}
\end{theorem}

\cref{thm:sd:sr,thm:hsd:mask,thm:hsd:converg} together provide us a comprehensive picture on Hermite subdivision schemes.
In \cref{sec:converg}, we shall prove \cref{thm:hsd:converg} and discuss factorization of general matrix-valued masks and Hermite masks through the approach of the normal form of matrix-valued masks.

The structure of this paper is as follows. In \cref{sec:sd}, we shall study how a general vector subdivision operator acts on vector polynomials. This provides us the basics for studying Hermite subdivision schemes and for proving \cref{thm:sd:sr}.
In \cref{sec:hsd}, we shall study
some necessary conditions for a convergent Hermite subdivision scheme.
Then we shall prove \cref{thm:hsd:mask} in \cref{sec:hsd}.
Finally, in \cref{sec:converg} we shall
first recall the normal form (also called the canonical form) of a matrix mask introduced and studied in \cite{han03,han09,hanbook,hm03}.
The notion of the normal form of a matrix-valued mask plays a central role for greatly facilitating the study of vector cascade algorithms and refinable vector functions.
For example, as a direct consequence of a normal form of a matrix-valued mask, we can straightforwardly obtain a
factorization of matrix-valued masks for Hermite and vector subdivision schemes.
Finally, we shall prove \cref{thm:hsd:converg} on the convergence of Hermite subdivision schemes.

\section{Action of Vector Subdivision Operators on Vector Polynomials}
\label{sec:sd}

To study various properties of Hermite and vector subdivision schemes, it is crucial to understand how a vector subdivision operator acts on vector polynomial spaces.
For scalar subdivision schemes, this problem is well studied in \cite{dl02,han03,han13,hanbook} and references therein. In this section we shall investigate vector subdivision operators acting on vector polynomials.
As we shall explain later in this section, it is much more complicated to study vector subdivision operators acting on vector polynomials than their scalar counterparts in \cite{han03,han13,hanbook}. In particular, we shall prove \cref{thm:sd:sr} in this section.

The following result is known in
\cite[Proposition~2.1]{han09} and \cite[Lemma~3.1]{han13}. For the convenience of the reader, we provide a proof here.

\begin{lemma}\label{lem:conv}
Let $u\in \lp{0}$ be a finitely supported sequence and $\pp$ be a polynomial of degree $m$. Then $\pp*u$ is a polynomial sequence satisfying $\deg(\pp*u)\le m$,
\be \label{polyconv}
(\pp*u)(x):=\sum_{k\in \Z} \pp(x-k)u(k)=\sum_{j=0}^\infty \frac{(-i)^j}{j!} \pp^{(j)}(x) \wh{u}^{(j)}(0),\qquad x\in \Z
\ee
and $[\pp*u]^{(k)}(\cdot-\tau)=(\pp^{(k)}(\cdot-\tau))*u$ for all $k\in \NN$ and $\tau\in \R$. Moreover, for $v\in \lp{0}$, the identity $\pp*u=\pp*v$ holds if and only if
$\wh{u}(\xi)=\wh{v}(\xi)+\bo(|\xi|^{m+1})$ as $\xi \to 0$.
\end{lemma}

\bp By the Taylor expansion $\pp(x-k)=\sum_{j=0}^\infty \pp^{(j)}(x)\frac{(-k)^j}{j!}=\sum_{j=0}^m \pp^{(j)}(x)\frac{(-k)^j}{j!}$, we have
\[
(\pp*u)(x)=\sum_{k\in \Z} \pp(x-k)u(k)=
\sum_{k\in \Z} \sum_{j=0}^\infty
\pp^{(j)}(x) u(k) \frac{(-k)^j}{j!}
=\sum_{j=0}^\infty \frac{1}{j!} \pp^{(j)}(x)
\sum_{k\in \Z} u(k) (-k)^j.
\]
Since $\wh{u}(\xi)=\sum_{k\in \Z} u(k) e^{-ik\xi}$, we have $\wh{u}^{(j)}(0)=\sum_{k\in \Z} u(k) (-ik)^j=i^j \sum_{k\in \Z} u(k) (-k)^j$. Hence, we conclude from the above identity that
\[
(\pp*u)(x)=\sum_{j=0}^\infty \frac{1}{j!} \pp^{(j)}(x) (-i)^j \wh{u}^{(j)}(0)=
\sum_{j=0}^\infty \frac{(-i)^j}{j!} \pp^{(j)}(x) \wh{u}^{(j)}(0).
\]
This proves \eqref{polyconv} which also holds for all $x\in \R$. For $k\in \NN$, using \eqref{polyconv} with $x\in \R$, we have
\[
[\pp*u]^{(k)}(x-\tau)=\sum_{j=0}^\infty \frac{(-i)^j}{j!} \pp^{(j+k)}(x-\tau) \wh{u}^{(j)}(0)=
\sum_{j=0}^\infty \frac{(-i)^j}{j!} [\pp(\cdot-\tau)]^{(j+k)}(x) \wh{u}^{(j)}(0)
=((\pp^{(k)}(\cdot-\tau))*u)(x).
\]
This proves $[\pp*u]^{(k)}(\cdot-\tau)=(\pp^{(k)}(\cdot-\tau))*u$.
Since $\pp,\pp',\ldots,\pp^{(m)}$ are linearly independent by $\deg(\pp)=m$,
we deduce from \eqref{polyconv} that $\pp*(u-v)=0$ if and only if $\wh{u-v}^{(j)}(0)=0$ for all $j=0,\ldots,m$.
\ep

For $m\in \NN:=\N\cup\{0\}$ and $v\in \lrs{0}{1}{r}$, recall that $\PR_{m,v}:=\{\pp*v \setsp \pp\in \PL_m\}$ in \eqref{Pvm}.
By \cref{lem:conv}, the space $\PR_{m,v}$ is invariant under derivatives and translation, that is,
\be \label{Pmv:translation}
\vec{\pp}^{(j)}(\cdot-\tau)\in \PR_{m,v}\qquad \forall\, \vec{\pp}\in \PR_{m,v}\;\; \mbox{and}\;\;  \tau \in \R,  j\in \NN.
\ee

The following result facilitates the study of subdivision operators acting on vector polynomials.

\begin{lemma}\label{lem:vecp}
Let $\vec{\pp}\in (\PL_m)^{1\times r}$ and $v\in \lrs{0}{1}{r}$. Then $\vec{\pp}=\frac{(\cdot)^m}{m!}*v$
if and only if
\be \label{poly:v}
\wh{v}^{(j)}(0)=j! i^j \vec{\pp}^{(m-j)}(0),\qquad j=0,\ldots,m.
\ee
Moreover, if $\vec{\pp}=\frac{(\cdot)^m}{m!}*v$, then
\be \label{pvm}
\PR_{\vec{\pp}^{(j)}}:=\mspan\{\vec{\pp}^{(j)},\vec{\pp}^{(j+1)},\ldots, \vec{\pp}^{(m)}\}=\PR_{m-j,v},\qquad j=0,\ldots,m.
\ee
\end{lemma}

\bp Note that $[\frac{(\cdot)^m}{m!}]^{(j)}=\frac{(\cdot)^{m-j}}{(m-j)!}$ for all $j=0,\ldots,m$.
Considering the Taylor expansion of $\vec{\pp}$ at the point $0$,
we have $\vec{\pp}(x)=\sum_{j=0}^m \vec{\pp}^{(j)}(0) \frac{x^j}{j!}$.
By \eqref{polyconv} in \cref{lem:conv},
\[
\left(\frac{(\cdot)^m}{m!}*v\right)(x)
=\sum_{j=0}^m \frac{(-i)^j}{j!} \frac{x^{m-j}}{(m-j)!} \wh{v}^{(j)}(0)
=\sum_{j=0}^m \frac{(-i)^{m-j}}{(m-j)!}
\frac{x^j}{j!} \wh{v}^{(m-j)}(0).
\]
Since $\frac{(\cdot)^j}{j!},j=0,\ldots,m$ are linearly independent and $\vec{\pp}(x)=\sum_{j=0}^m \vec{\pp}^{(j)}(0)\frac{x^j}{j!}$, we conclude from the above identity that $\vec{\pp}=\frac{(\cdot)^m}{m!}*v$ if and only if \eqref{poly:v} holds.

By \eqref{lem:conv} and $\vec{\pp}=\frac{(\cdot)^m}{m!}*v$, we have
$\vec{\pp}^{(j)}=[\frac{(\cdot)^m}{m!}*v]^{(j)}=
\frac{(\cdot)^{m-j}}{(m-j)!}*v\in \PR_{m-j,v}$ for all $j\in \NN$.
This proves $\PR_{\vec{\pp}^{(j)}}
\subseteq \PR_{m-j,v}$.
Conversely, any $\pq\in \PL_{m-j}$ has the Taylor expansion $\pq(x)=\sum_{k=0}^{m-j} \pq^{(k)}(0) \frac{x^k}{k!}$. Hence,
\[
(\pq*v)(x)=\sum_{k=0}^{m-j} \pq^{(k)}(0)\left( \frac{(\cdot)^k}{k!}*v\right)(x)=
\sum_{k=0}^{m-j} \pq^{(k)}(0) \vec{\pp}^{(m-k)}(x)
\]
which belongs to $\PR_{\vec{\pp}^{(j)}}$. Hence, $\PR_{m-j,v}\subseteq \PR_{\vec{\pp}^{(j)}}$.
This proves \eqref{pvm}.
\ep

For $a\in \lrs{0}{r}{r}$ and $\gamma\in\Z$, recall that $a^{[\gamma]}:=\{a(\gamma+2k) \}_{k\in \Z}$ is its \emph{$\gamma$-coset sequence}. Note that $\wh{a^{[\gamma]}}(\xi)=\sum_{k\in \Z} a(\gamma+2k)e^{-ik\xi}$ and $\wh{a}(\xi)=\wh{a^{[0]}}(2\xi)+\wh{a^{[1]}}(2\xi)e^{-i\xi}$.

Generalizing \cite[Theorem~3.4]{han13} and \cite[Theorem~1.2.4]{hanbook}
from scalar subdivision operators to vector subdivision operators,
we have the following result about a vector subdivision operator acting on a given vector polynomial.

\begin{theorem}\label{thm:sd:poly}
Let $\vec{\pp}\in (\PL_m)^{1\times r}$ and $v\in \lrs{0}{1}{r}$ such that
$\vec{\pp}=\frac{(\cdot)^m}{m!}*v$.
Define $u_a\in \lrs{0}{1}{r}$ by $\wh{u_a}(\xi):=\wh{v}(2\xi)\wh{a}(\xi)$.
Then the following statements are equivalent to each other:
\begin{enumerate}
\item[(1)] $\sd_a \vec{\pp}\in (\PL_m)^{1\times r}$, that is, $\sd_a \vec{\pp}$ is a vector polynomial sequence.
\item[(2)] $\wh{v}(\xi)\wh{a^{[1]}}(\xi) e^{-i\xi/2}=\wh{v}(\xi)\wh{a^{[0]}}(\xi)+\bo(|\xi|^{m+1})$ as $\xi \to 0$.
\item[(3)] $\wh{v}(2\xi) \wh{a}(\xi+\pi)=\bo(|\xi|^{m+1})$ as $\xi \to 0$, i.e., $\wh{u_a}(\xi+\pi)=\bo(|\xi|^{m+1})$ as $\xi \to 0$.
\item[(4)]
$\sd_a \PR_{m,v}=\PR_{m,u_a}$. Note that $\PR_{m,v}=\PR_{\vec{\pp}}:=\mspan\{\vec{\pp},\vec{\pp}',\ldots,\vec{\pp}^{(m)}\}$.
\item[(5)] $\sd_a (\vec{\pp}^{(j)})\in (\PL_m)^{1\times r}$ for all $j\in \NN$, i.e., $\sd_a \PR_{\vec{\pp}}\subseteq (\PL_m)^{1\times r}$.
\end{enumerate}
Moreover, any of the above items (1)--(5) implies
\be \label{sd:poly}
\sd_a (\pq*v) =(\pq(2^{-1}\cdot))*u_a,
\quad \forall\, \pq\in \PL_m,
\quad \mbox{or equivalently},\quad
\sd_a \vec{\pq}=(\vec{\pq}(2^{-1}\cdot))*a,\quad
\forall\, \vec{\pq}\in \PR_{m,v},
\ee
and
\be \label{sd:poly:2}
\sd_a (\vec{\pq}^{(j)})=2^j [\sd_a \vec{\pq}]^{(j)},\qquad \forall\, \vec{\pq}\in \PR_{m,v}, j\in \NN.
\ee
\end{theorem}

\bp
Define $\pp_m(x):=\frac{x^m}{m!}$. By assumption, we have $\vec{\pp}=\pp_m*v$ and hence $\PR_{m,v}=\PR_{\vec{\pp}}$ by \cref{lem:vecp}.
For $\gamma\in \Z$,
note that $u_a^{[\gamma]}=v*a^{[\gamma]}$.
By \cref{lem:conv} and the definition of $\sd_a$ in \eqref{sd}, we have
\begin{align*}
(\sd_a \vec{\pp})(\gamma+2 \cdot)
&=2\sum_{k\in \Z} \vec{\pp}(k) a(\gamma+2\cdot-2k)=
2\sum_{k\in \Z} \vec{\pp}(k) a^{[\gamma]}(\cdot-k)
=
2\vec{\pp}*a^{[\gamma]}\\
&=2(\pp_m*v)*a^{[\gamma]}
=2\pp_m*u_a^{[\gamma]}=
2\sum_{j=0}^\infty \frac{(-i)^j}{j!} \pp_m^{(j)}(\cdot) \wh{u_a^{[\gamma]}}^{(j)}(0).
\end{align*}
Hence, we deduce from the above identity that $\sd_a \vec{\pp}(x)=\vec{\theta}(x)$ for all $x\in (\gamma+2\Z)$, where
\be \label{sd:poly:eq1}
\vec{\theta}(x):=2\sum_{j=0}^\infty \frac{(-i)^j}{j!} \pp_m^{(j)}(2^{-1}x-2^{-1}\gamma)\wh{u_a^{[\gamma]}}^{(j)}(0).
\ee
Note that $\sd_a \vec{\pp}\in (\PL_m)^{1\times r}$ if and only if $\vec{\theta}$ is a vector polynomial independent of the choice of $\gamma\in \Z$.
We now calculate $\vec{\theta}$. Using the Taylor expansion of $\pp^{(j)}(2^{-1}x-2^{-1}\gamma)$ at $2^{-1}x$, we have
$\pp^{(j)}(2^{-1} x-2^{-1} \gamma)=\sum_{n=0}^\infty \frac{1}{n!} \pp^{(j+n)}(2^{-1}x) (-2^{-1}\gamma)^n=\sum_{n=0}^m \frac{1}{n!} \pp^{(j+n)}(2^{-1}x) (-2^{-1}\gamma)^n$. Using substitution $k=j+n$, we have
\begin{align*}
\sum_{j=0}^\infty \frac{(-i)^j}{j!} \pp_m^{(j)}(2^{-1}x-2^{-1}\gamma) \wh{u_a^{[\gamma]}}^{(j)}(0)
&=
\sum_{j=0}^\infty \sum_{n=0}^\infty \frac{(-i)^j}{j! n!}\pp_m^{(j+n)}(2^{-1}x) (-2^{-1}\gamma)^n \wh{u_a^{[\gamma]}}^{(j)}(0)\\
&=\sum_{k=0}^\infty
\sum_{j=0}^k \frac{(-i)^j}{j!(k-j)!}
\pp_m^{(k)}(2^{-1}x)(-2^{-1}\gamma)^{k-j} \wh{u_a^{[\gamma]}}^{(j)}(0)\\
&=\sum_{k=0}^\infty \pp_m^{(k)}(2^{-1}x) \sum_{j=0}^k \frac{(-i)^j}{j!(k-j)!} (-2^{-1}\gamma)^{k-j} \wh{u_a^{[\gamma]}}^{(j)}(0).
\end{align*}
On the other hand, using Leibniz differentiation formula, we have
\begin{align*}
\Big[\wh{u_a^{[\gamma]}}(\xi) e^{-i\gamma\xi/2}\Big]^{(k)}(0)
&=\sum_{j=0}^k \frac{k!}{j!(k-j)!}
[e^{-i\gamma\xi/2}]^{(k-j)}(0)\wh{u_a^{[\gamma]}}^{(j)}(0)\\
&=\sum_{j=0}^k \frac{k!}{j!(k-j)!}
(-i2^{-1}\gamma)^{k-j}\wh{u_a^{[\gamma]}}^{(j)}(0)
=\sum_{j=0}^k \frac{k! i^{k-j}}{j!(k-j)!}
(-2^{-1}\gamma)^{k-j}\wh{u_a^{[\gamma]}}^{(j)}(0).
\end{align*}
Combining the above two identities, we conclude that
\[
\sum_{j=0}^\infty \frac{(-i)^j}{j!} \pp_m^{(j)}(2^{-1}x-2^{-1}\gamma) \wh{u_a^{[\gamma]}}^{(j)}(0)
=\sum_{k=0}^\infty \pp_m^{(k)}(2^{-1}x) \frac{(-i)^k}{k!} \Big[\wh{u_a^{[\gamma]}}(\xi) e^{-i\gamma\xi/2}\Big]^{(k)}(0).
\]
In summary, the identity in \eqref{sd:poly:eq1} becomes
\be \label{sd:poly:eq2}
\vec{\theta}(x)=2\sum_{j=0}^m \frac{(-i)^j}{j!}  \pp_m^{(j)}(2^{-1}x) \Big[\wh{u_a^{[\gamma]}}(\xi) e^{-i\gamma\xi/2}\Big]^{(j)}(0).
\ee
Note that $\pp_m,\pp_m',\ldots,\pp_m^{(m)}$ are linearly independent by $\pp_m(x)=\frac{x^{m}}{m!}$. Therefore, it follows from \eqref{sd:poly:eq2} that $\vec{\theta}$ is a vector polynomial independent of $\gamma\in \Z$ if and only if $[\wh{u_a^{[\gamma]}}(\xi)e^{-i\gamma\xi/2}]^{(j)}(0)$ is independent of $\gamma\in \{0,1\}$ for all $j=0,\ldots,m$. Since $\wh{u_a^{[\gamma]}}(\xi)
=\wh{v}(\xi)\wh{a^{[\gamma]}}(\xi)$ and $\sd_a \vec{\pp}=\vec{\theta}$,
this proves (1)$\iff$(2).

We now prove that item (2) implies \eqref{sd:poly} and \eqref{sd:poly:2}. We first prove \eqref{sd:poly} with $\pq(x)=\pp_m(x)=\frac{x^m}{m!}$. By what has been proved, we have $\sd_a \vec{\pp}=\vec{\theta}$, where $\vec{\theta}$ in \eqref{sd:poly:eq2} is a vector polynomial independent of $\gamma$. Averaging \eqref{sd:poly:eq2} for $\gamma\in \{0,1\}$ and noting that $\wh{u_a}(\xi/2)=\wh{u_a^{[0]}}(\xi)+\wh{u_a^{[1]}}(\xi)e^{-i\xi/2}$, we obtain
\begin{align*}
(\sd_a \vec{\pp})(x)
&=\sum_{j=0}^m \frac{(-i)^j}{j!}
\pp_m^{(j)}(2^{-1}x) \sum_{\gamma=0}^1
\Big[\wh{u_a^{[\gamma]}}(\xi) e^{-i\gamma\xi/2}\Big]^{(j)}(0)
=\sum_{j=0}^m \frac{(-i)^j}{j!} \pp_m^{(j)}(2^{-1}x)  [\wh{u_a}(\cdot/2)]^{(j)}(0)\\
&=\sum_{j=0}^m \frac{(-i)^j}{j!} [\pp_m(2^{-1}\cdot)]^{(j)}(x) \wh{u_a}^{(j)}(0)=(\pp_m(2^{-1}\cdot))*u_a.
\end{align*}
Since $\vec{\pp}=\frac{(\cdot)^m}{m!}*v=\pp_m*v$, we proved that item (2) implies
\be \label{sd:poly:eq3}
\sd_a (\pp_m*v)=(\pp_m(2^{-1}\cdot))*u_a.
\ee
On the other hand, item (2) trivially implies
\[
\wh{v}(\xi) \wh{a^{[1]}}(\xi)e^{-i\xi/2}=
\wh{v}(\xi)\wh{a^{[0]}}(\xi)+\bo(|\xi|^{j+1}),\quad
\xi \to 0\quad \mbox{for all}\quad j=0,\ldots,m.
\]
That is, item (2) holds with $m$ being replaced by $j$ for all $j=0,\ldots,m$.
By what has been proved, \eqref{sd:poly:eq3} must hold with $\pp_m$ being replaced by $\pp_j$ with $\pp_j(x):=\frac{x^j}{j!}$ for all $j=0,\ldots,m$. Since $\mspan\{\pp_0,\ldots,\pp_m\}=\PL_m$, this proves the first identity in \eqref{sd:poly}.

Since $\wh{u_a}(\xi)=\wh{v}(2\xi)\wh{a}(\xi)$, by Leibniz differentiation formula, we have
\be \label{ua}
\wh{u_a}^{(j)}(0)=\sum_{k=0}^j \frac{j!}{k!(j-k)!} 2^{j-k}\wh{v}^{(j-k)}(0)\wh{a}^{(k)}(0).
\ee
To prove the second identity in \eqref{sd:poly}, for $\vec{\pq}\in \PR_{m,v}$, there must exist $\pq\in \PL_m$ such that $\vec{\pq}=\pq*v$. Therefore, by the first identity of \eqref{sd:poly} and \eqref{polyconv}, using \eqref{ua}, we deduce that
\begin{align*}
\sd_a\vec{\pq}
&=\sd_a(\pq*v)=(\pq(2^{-1}\cdot))*u_a
=\sum_{j=0}^\infty \frac{(-i)^j}{j!} [\pq(2^{-1}\cdot)]^{(j)}(\cdot) \wh{u_a}^{(j)}(0)\\
&=\sum_{j=0}^\infty\sum_{k=0}^j  \frac{(-i)^j}{k!(j-k)!} 2^{-k} \pq^{(j)}(2^{-1}\cdot) \wh{v}^{(j-k)}(0)\wh{a}^{(k)}(0)\\
&=\sum_{k=0}^\infty
\frac{(-i)^k}{k!}2^{-k}
\left(\sum_{j=k}^\infty \frac{(-i)^{j-k}}{(j-k)!}[\pq^{(k)}]^{(j-k)}(2^{-1}\cdot)
\wh{v}^{(j-k)}(0)\right) \wh{a}^{(k)}(0)\\
&=\sum_{k=0}^\infty \frac{(-i)^k}{k!}
2^{-k}(\pq^{(k)}*v)(2^{-1}\cdot) \wh{a}^{(k)}(0)
=\sum_{k=0}^\infty \frac{(-i)^k}{k!}
[\vec{\pq}(2^{-1}\cdot)]^{(k)} \wh{a}^{(k)}(0)
=(\vec{\pq}(2^{-1}\cdot))*a,
\end{align*}
where we used $[\vec{\pq}(2^{-1}\cdot)]^{(k)}=
2^{-k} \vec{\pq}^{(k)}(2^{-1}\cdot)
=2^{-k} (\pq^{(k)}*v)(2^{-1}\cdot)$.
This proves \eqref{sd:poly}. To prove that item (2) implies \eqref{sd:poly:2}, by proved \eqref{sd:poly}, we have
\[
2^j [\sd_a\vec{\pq}]^{(j)}=
2^j [\vec{\pq}(2^{-1}\cdot)*a]^{(j)}
=(\vec{\pq}^{(j)}(2^{-1}\cdot))*a
=\sd_a (\vec{\pq}^{(j)}).
\]
This proves \eqref{sd:poly:2}.
Hence, item (2) implies both \eqref{sd:poly} and \eqref{sd:poly:2}.

Note that $\wh{a}(\xi+\pi)=\wh{a^{[0]}}(2\xi)-\wh{a^{[1]}}(2\xi)e^{-i\xi}$. Item (3) is equivalent to $\wh{v}(2\xi)\wh{a^{[1]}}(2\xi)e^{-i\xi}=\wh{v}(2\xi)
\wh{a^{[0]}}(2\xi)+\bo(|\xi|^{m+1})$ as $\xi \to 0$, which is obviously equivalent to item (2). This proves (2)$\iff$(3).

Because we already proved (3)$\iff$(2)$\iff$(1),
(3)$\imply$(4) follows directly from the first identity in \eqref{sd:poly}.
Since $\vec{\pp}^{(j)}=\frac{(\cdot)^{m-j}}{(m-j)!}*v\in \PR_{m,v}$, it is trivial that (4)$\imply$(5). Item (5)$\imply$(1) is obvious. This proves all the claims.
\ep

Let us now look at the differences between  $r=1$ and $r>1$.
Suppose that $\vec{\pp}$ is a vector polynomial with $m:=\deg(\vec{\pp})$ such that $\sd_a\vec{\pp}\in (\PL_m)^{1\times r}$. We now consider two cases.

Case 1: $r=1$, i.e., we have a scalar subdivision scheme.
Then it is trivial to observe that $\mspan\{\vec{\pp},\vec{\pp}',\ldots,\vec{\pp}^{(m)}\}=\PL_m$ and $\PR_{m,u_a}\subseteq\PL_m$. Consequently, $\sd_a^n \PL_m\subseteq \PL_m$ for all $n\in \N$.

Case 2: $r>1$, i.e., we have a vector subdivision scheme. Then $\PR_{\vec{\pp}}:=
\mspan\{\vec{\pp},\vec{\pp}',\ldots,\vec{\pp}^{(m)}\}
\subsetneq (\PL_m)^{1\times r}$. That is, $\PR_{\vec{\pp}}$ must be a proper subspace of $(\PL_m)^{1\times r}$ and it is not necessary that $\PR_{m,u_a}\subseteq \PR_{\vec{\pp}}$. Consequently, even if $\sd_a \vec{\pp}\in (\PL_m)^{1\times r}$, we no longer necessarily have $\sd_a^2 \vec{\pp} \in (\PL_m)^{1\times r}$.
Since a vector subdivision scheme is always iterative,
the polynomial reproduction property may lose during iteration. This makes things a lot complicated for vector subdivision schemes.

Next, let us study when $\sd_a^n \vec{\pp}\in (\PL_m)^{1\times r}$ for all $n=1,\ldots,N$ with $N\in \N\cup\{\infty\}$.

\begin{cor}\label{cor:sd}
Let $\vec{\pp}\in (\PL_m)^{1\times r}$ and $v\in \lrs{0}{1}{r}$ such that $\vec{\pp}=\frac{(\cdot)^m}{m!}*v$. Let $N\in \N\cup\{\infty\}$ and finitely supported masks $a_n\in \lrs{0}{r}{r}$ for $n=1,\ldots,N$.
Then the following are equivalent:
\begin{enumerate}
\item[(1)] $\sd_{a_n}\cdots \sd_{a_1} \vec{\pp}\in (\PL_m)^{1\times r}$ for all $n=1,\ldots,N$.
\item[(2)] $\wh{v_n}(\xi+\pi)=\bo(|\xi|^{m+1})$ as $\xi \to 0$ for all $n=1,\ldots, N$, where
\[
\wh{v_n}(\xi):=\wh{v}(2^{n}\xi)\wh{a_1}(2^{n-1}\xi)
    \wh{a_2}(2^{n-2}\xi)
    \cdots\wh{a_{n-1}}(2\xi)\wh{a_n}(\xi),\qquad n\in \N.
\]

\item[(3)] $\sd_{a_n}\cdots \sd_{a_1} \PR_{m,v}=\PR_{m,v_n}$ for all $n=1,\ldots,N$.
\end{enumerate}
\end{cor}

\bp Note that $\wh{v_n}(\xi)=\wh{v_{n-1}}(2\xi)\wh{a_{n}}(\xi)$. Now all the claims follow directly from \cref{thm:sd:poly}.
\ep

Let $a_n\in \lrs{0}{r}{r}$ for $n\in \N$. For $N\in \N$, we define
\[
\PR_{N}:=\{\vec{\pp}\in (\PL_m)^{1\times r} \setsp \sd_{a_n}\cdots \sd_{a_1} \vec{\pp}\in (\PL_m)^{1\times r}\quad \mbox{for all}\quad n=1,\ldots,N\}
\]
and $\PR_{\infty}:=\cap_{N=1}^\infty \PR_{N}$. Note that $\PR_\infty\subseteq \PR_{N+1}\subseteq \PR_{N}\subseteq (\PL_m)^{1\times r}$ for all $N\in \N$.
Using \cref{lem:vecp}, we observe that
the linear space $\PR_{N}$ consists of all elements $\sum_{j=0}^m \frac{(-i)^j}{j!} \frac{x^{m-j}}{(m-j)!} \wh{v}^{(j)}(0)$, where $\{\wh{v}^{(j)}(0)\}_{j=0}^{m}$ is a solution to the linear systems induced by item (2) of \cref{cor:sd}.

If $a_n=a$ for all $n\in \N$,
since $\{\dim(\PR_N)\}_{N=1}^\infty$ is a decreasing sequence, then there must exist $N\in \N$ such that
$\PR_{N+1}=\PR_{N}$ and hence $\PR_\infty=\PR_N$.
Since $\PR_\infty$ is closed under derivatives and translation (see \eqref{Pmv:translation}), if $\PR_\infty\ne \{0\}$, then $\PR_\infty$ must contain a nonzero constant vector sequence.
For a convergent subdivision scheme, it is natural to require that $\sd_a c=c$ for some nonzero constant vector sequence $c\in \PR_\infty$, that is, the constant vector sequence $c$ must be preserved by the subdivision operator $\sd_a$.

We are now ready to prove the main result stated in \cref{thm:sd:sr} on a vector subdivision operator acting on vector polynomials.

\bp[Proof of \cref{thm:sd:sr}] By \cref{lem:vecp}, we have $\PR_{\vec{\pp}}:=\mspan\{\vec{\pp},\vec{\pp}',\ldots,\vec{\pp}^{(m)}\}=\PR_{m,v}$ and $\vec{\pp}^{(m)}=1*v=\wh{v}(0)$. Therefore, it is obvious that (1)$\imply$(2).

(2)$\imply$(3). Since $\sd_a\vec{\pp}\in \PR_{\vec{\pp}}\subseteq (\PL_m)^{1\times r}$, by \cref{thm:sd:poly}, the second identity in \eqref{sr:1} and the identities in \eqref{sd:poly} must hold. By $\vec{\pp}_m=1*v=\wh{v}(0)$ and the first identity in \eqref{sd:poly}, we have $\sd_a (\vec{\pp}^{(m)})=1*u_a=\wh{u_a}(0)=\wh{v}(0)\wh{a}(0)$, where $\wh{u_a}(\xi):=\wh{v}(2\xi)\wh{a}(\xi)$.
By our assumption $\sd_a (\vec{\pp}^{(m)})=\vec{\pp}^{(m)}=\wh{v}(0)$, we conclude that $\wh{v}(0)\wh{a}(0)=\sd_a (\vec{\pp}^{(m)})=\wh{v}(0)$.
Using \eqref{sd:poly} and \cref{lem:conv}, we have
\be \label{sd:pp}
\sd_a \vec{\pp}=
\sd_a(\pq_m*v)=(\pq_m(2^{-1}\cdot))*u_a=
2^{-m} \sum_{n=0}^m \frac{(-i)^{n}}{n!}\pq_{m-n}(\cdot)\wh{u_a}^{(n)}(0),
\ee
where $\pq_{m-n}(x):=\frac{x^{m-n}}{(m-n)!}$ for $n=0,\ldots,m$.
By our assumption $\sd_a \vec{\pp}\in \mspan\{\vec{\pp},\vec{\pp}',\ldots,\vec{\pp}^{(m)}\}$, we have $\sd_a\vec{\pp}=\sum_{k=0}^m d_k \vec{\pp}^{(k)}$ for some $d_0,\ldots,d_m\in \C$. By \cref{lem:vecp}, $\vec{\pp}^{(k)}(x)=\sum_{j=0}^{m-k} \frac{(-i)^j}{j!} \frac{x^{m-k-j}}{(m-k-j)!} \wh{v}^{(j)}(0)$. Therefore, using substitution $n=k+j$, we have
\[
\sd_a\vec{\pp}(x)=\sum_{k=0}^m d_k \vec{\pp}^{(k)}(x)=
\sum_{k=0}^m
\sum_{j=0}^{m-k} d_k \frac{(-i)^j}{j!} \frac{x^{m-k-j}}{(m-j-k)!} \wh{v}^{(j)}(0)
=\sum_{n=0}^m \pq_{m-n}(x)\sum_{j=0}^{n} d_{n-j} \frac{(-i)^j}{j!} \wh{v}^{(j)}(0).
\]
Comparing \eqref{sd:pp} with the above identity and noting that $\{\pq_{m-n}\}_{n=0}^m$ are linearly independent, we must have
$2^{-m}\frac{(-i)^{n}}{n!}\wh{u_a}^{(n)}(0)
=\sum_{j=0}^{n} d_{n-j} \frac{(-i)^j}{j!} \wh{v}^{(j)}(0)$
for all $n=0,\ldots,m$. That is,
\be \label{va:v}
\wh{u_a}^{(n)}(0)=\sum_{j=0}^n i^{n-j} 2^{m} (n-j)!d_{n-j} \frac{n!}{j!(n-j)!} \wh{v}^{(j)}(0),\qquad n=0,\ldots,m.
\ee
Take $c\in \lp{0}$ such that $\wh{c}^{(j)}(0)=2^{m} i^j j! d_j$ for all $j=0,\ldots,m$. Then
$i^{n-j}2^m (n-j)!d_{n-j}=\wh{c}^{(n-j)}(0)$.
By Leibniz differentiation formula, \eqref{va:v} is equivalent to $\wh{u_a}(\xi)=\wh{c}(\xi)\wh{v}(\xi)+\bo(|\xi|^{m+1})$ as $\xi\to 0$. Because $\wh{u_a}(\xi)=\wh{v}(2\xi)\wh{a}(\xi)$, this proves the first identity in \eqref{sr:1}.
Taking $n=0$ in \eqref{va:v}, we have $\wh{u_a}(0)=2^{m}d_0\wh{v}(0)$. Since we proved $\wh{u_a}(0)=\wh{v}(0)\wh{a}(0)=\wh{v}(0)$, we conclude from $2^{m} d_0\wh{v}(0)=\wh{u_a}(0)=\wh{v}(0)$ and $\wh{v}(0)\ne 0$ that $d_0=2^{-m}$. Consequently, $\wh{c}(0)=2^{m}d_0=1$.
This proves (2)$\imply$(3).

(3)$\imply$(4). Since $c\in \lp{0}$ with $\wh{c}(0)=1$, the function $\wh{\varphi}(\xi):=\prod_{j=1}^\infty \wh{c}(2^{-j}\xi)$ is a well-defined analytic function and satisfies $\wh{\varphi}(2\xi)=\wh{c}(\xi)\wh{\varphi}(\xi)$ with $\wh{\varphi}(0)=1$. Note that $\wh{c}(\xi)=\wh{\varphi}(2\xi)/\wh{\varphi}(\xi)$.
Take a sequence $d\in \lp{0}$ satisfying $\wh{d}(\xi)=1/\wh{\varphi}(\xi)+\bo(|\xi|^{m+1})$ as $\xi \to 0$. Note that $\wh{d}(0)=1/\wh{\varphi}(0)=1$.
By $\wh{\vgu_a}(\xi):=\wh{d}(\xi)\wh{v}(\xi)=
\wh{v}(\xi)/\wh{\varphi}(\xi)+\bo(|\xi|^{m+1})$ as $\xi\to 0$, we can easily see that
\eqref{sr:1} implies
\eqref{sr}.

(4)$\imply$(5). Since $\wh{\vgu_a}(\xi)=\wh{d}(\xi)\wh{v}(\xi)$ for some $d\in \lp{0}$ with $\wh{d}(0)=1$, we have $\PR_{m,v}=\PR_{m,\vgu_a}$.
By item (4), all the claims in Theorem~\ref{thm:sd:poly} hold with $\vec{\pp}=\vec{\pp}_m:=\frac{(\cdot)^m}{m!}*\vgu_a$.
Note that $\vec{\pp}_j:=\frac{(\cdot)^j}{j!}*\vgu_a
=\vec{\pp}_m^{(m-j)}$ is obviously true.
We conclude from the first identity in \eqref{sd:poly} that
\[
\sd_a \vec{\pp}_j=\sd_a \left (\frac{(\cdot)^j}{j!}*\vgu_a\right)
=
\left(\frac{(2^{-1}\cdot)^j}{j!}\right)*u_a=
2^{-j} \left(\frac{(\cdot)^j}{j!}\right)*u_a,
\]
where $\wh{u_a}(\xi):=\wh{\vgu_a}(2\xi)\wh{a}(\xi)$.
By the first identity in \eqref{sr}, we have $\wh{u_a}(\xi)=\wh{\vgu_a}(\xi)+\bo(|\xi|^{m+1})$ as $\xi \to 0$. Consequently, we conclude from the above identity and \cref{lem:conv} that for all $j=0,\ldots,m$,
\[
\sd_a \vec{\pp}_j=2^{-j} \left(\frac{(\cdot)^j}{j!}\right)*u_a
=2^{-j} \left(\frac{(\cdot)^j}{j!}\right)*\vgu_a
=2^{-j}\vec{\pp}_j.
\]
This proves (4)$\imply$(5).

(5)$\imply$(6) is trivial since $\vec{\pp}_m$ in \eqref{sd:eigpoly} belongs to $\PR_{m,\vgu_a}$ and $\deg(\vec{\pp}_m)=m$ by $\wh{v}(0)\ne 0$. Note that $\PR_{m,\vgu_a}=\PR_{m,v}$ by Lemma~\ref{lem:vecp}.

(6)$\imply$(1). Since $\sd_a \vec{\pp}_m=2^{-m} \vec{\pp}_m\in (\PL_m)^{1\times r}$, by \cref{thm:sd:poly}, we have $\sd_a(\vec{\pp}_m^{(j)})=2^j [\sd_a\vec{\pp}_m]^{(j)}=2^{j-m} \vec{\pp}_m^{(j)}$. Since $\vec{\pp}_m\in \PR_{m,v}$ and $\PR_{m,v}$ is closed under derivatives, we have $\sd_a(\vec{\pp}_m^{(j)})=2^{j-m}\vec{\pp}_m^{(j)}\in \PR_{m,v}$.
In particular, $\sd_a (\vec{\pp}_m^{(m)})=\vec{\pp}_m^{(m)}$. Note that $\vec{\pp}_m^{(m)}=[\frac{(\cdot)^m}{m!}*v]^{(m)}=
1*v=\wh{v}(0)$. Hence, we proved $\sd_a (\wh{v}(0))=\wh{v}(0)$.
Since $\PR_{\vec{\pp}_m}=\mspan\{\vec{\pp}_m,\vec{\pp}_m',\ldots,\vec{\pp}_m^{(m)}\}$,
we have $\sd_a \PR_{\vec{\pp}_m}= \PR_{\vec{\pp}_m}\subseteq \PR_{m,v}$. Because $\deg(\vec{\pp}_m)=m$ and all the elements $\vec{\pp}_m,\vec{\pp}_m',\ldots,\vec{\pp}_m^{(m)}$ are linearly independent,  we deduce from $\PR_{\vec{\pp}_m}\subseteq \PR_{m,v}$ that $\PR_{\vec{\pp}_m}=\PR_{m,v}$.
Therefore, we proved $\sd_a \PR_{m,v}=\PR_{m,v}$.
This proves (6)$\imply$(1).

The identity in \eqref{sd:deriv} follows directly from
\eqref{sd:poly} and \eqref{sd:poly:2}.
\ep

As a special case of \cref{thm:sd:sr}, we have the following result.

\begin{cor}\label{cor:sr}
Let $a\in \lrs{0}{r}{r}$ and $m\in \NN$.
\begin{enumerate}
\item[(i)] If $\sd_a \vec{\pp}_m=2^{-m}\vec{\pp}_m$ for some $\vec{\pp}_m\in (\PL_m)^{1\times r}$ with $\deg(\vec{\pp}_m)=m$, then $a$ has order $m+1$ sum rules with respect to any sequence $\vgu_a\in \lrs{0}{1}{r}$ satisfying \eqref{sr} and $\wh{\vgu_a}^{(j)}(0)=j!i^j\vec{\pp}_m^{(m-j)}(0)$ for all $j=0,\ldots, m$. Note that $\wh{\vgu_a}(0)=\vec{\pp}_m^{(m)}(0)\ne 0$ due to $\deg(\vec{\pp}_m)=m$.
\item[(ii)] If $a$ has order $m+1$ sum rules in \eqref{sr} with respect to $\vgu_a\in \lrs{0}{1}{r}$ with $\wh{\vgu_a}(0)\ne 0$, then $\sd_a (\vec{\pp}_m^{(j)})=2^{j-m} \vec{\pp}_m^{(j)}$ and $\deg(\vec{\pp}_m^{(j)})=m-j$ for all $j=0,\ldots,m$, where $\vec{\pp}_m:=\frac{(\cdot)^m}{m!}*\vgu_a$.
\end{enumerate}
\end{cor}

\bp By the choice of $\vgu_a$ in item (i), we have $\vec{\pp}_m=\frac{(\cdot)^m}{m!}*\vgu_a$. Now items (i) and (ii) follow directly from \cref{thm:sd:sr}.
\ep


In terms of coset sequences, note that the definition of sum rules in \eqref{sr} can be equivalently expressed as (see \cite[Lemma~5.5.5]{hanbook}):
\be \label{sr:2}
\begin{split}
&\wh{\vgu_a}(\xi)\wh{a^{[0]}}(\xi)=2^{-1}\wh{\vgu_a}
(\xi/2)+\bo(|\xi|^{m+1}),\qquad \xi\to 0,\\
&\wh{\vgu_a}(\xi)\wh{a^{[1]}}(\xi)=2^{-1}e^{i\xi/2}
\wh{\vgu_a}
(\xi/2)+\bo(|\xi|^{m+1}), \qquad \xi\to 0,
\end{split}
\ee
where $a^{[\gamma]}(k):=a(\gamma+2k)$ for all $k\in \Z$ and $\wh{a^{[\gamma]}}(\xi):=\sum_{k\in \Z} a(\gamma+2k) e^{-ik\xi}$ for all $\gamma\in \Z$.

For the scalar case $r=1$, we must have $\wh{a}(0)=1$ in \eqref{sr} due to $\wh{\vgu_a}(0)\wh{a}(0)=\wh{\vgu_a}(0)$ and $\wh{\vgu_a}(0)\ne 0$. Moreover, we can pick $\vgu_a\in \lp{0}$ in \eqref{sr} through $\wh{\vgu_a}(\xi)=1/\wh{\varphi_a}(\xi)+\bo(|\xi|^{m+1})$ as $\xi\to 0$, where $\wh{\varphi_a}(\xi):=\prod_{j=1}^\infty \wh{a}(2^{-j}\xi)$. Then $\wh{\vgu_a}(0)=1$ and the first identity in \eqref{sr} is automatically satisfied. Hence, for the scalar case $r=1$, a mask $a$ has order $m+1$ sum rules in \eqref{sr} if and only if $\wh{a}(0)=1$
and $\wh{a}(\xi+\pi)=\bo(|\xi|^{m+1})$ as $\xi \to 0$, i.e., $\wh{a}(\xi)=2^{-1-m}(1+e^{-i\xi})^{m+1}
\wh{c}(\xi)$ for some $c\in \lp{0}$ with $\wh{c}(0)=1$. In addition, we have
$\sd_a (\vec{\pp}_a^{(j)})=2^{j-m}\vec{\pp}_a^{(j)}$ and $\deg(\vec{\pp}_a^{(j)})=m-j$ for all $j=0,\ldots,m$, where
\be \label{pa}
\pp_a(x):=\left(\frac{(\cdot)^m}{m!}*\vgu_a\right)(x)
=\sum_{j=1}^m \frac{(-i)^j}{j!} \frac{x^{m-j}}{(m-j)!} \wh{\vgu_a}^{(j)}(0)=
\sum_{j=1}^m \frac{(-i)^j}{j!} \frac{x^{m-j}}{(m-j)!} [1/\wh{\varphi_a}(\xi)]^{(j)}(0).
\ee
For the case $r>1$, things are a little bit more complicated and confusing. Let us provide an example here to demonstrate the differences between scalar and vector subdivision schemes.

\begin{example}{\rm
For $m\in \N$, the mask $a_m^B\in \lp{0}$ for the B-spline of order $m$ is given by $\wh{a_m^B}(\xi):=2^{-m}(1+e^{-i\xi})^m$.
Then $\wh{a^B_m}(0)=1$ and $a_m^B$ has order $m$ sum rules.
For $m,n\in \N$, we define a finitely supported matrix-valued mask $a\in \lrs{0}{2}{2}$ by
\[
\wh{a}(\xi):=\left[\begin{matrix} \wh{a_m^B}(\xi) &0\\
0 &2^{-m}\wh{a_n^B}(\xi)\end{matrix}\right].
\]
Define $\PR_\infty:=\cup_{m=0}^\infty \{ \pp\in (\PL_m)^{1\times 2} \setsp \sd_a^n \pp\in (\PL_m)^{1\times 2}\quad \forall\, n\in \N\}$. By the remark after \cref{cor:sr}, we have $\PR_\infty=\{(\pp,\pq)\setsp \pp\in \PL_{m-1},\pq\in \PL_{n-1}\}$ and $\sd_a \PR_\infty=\PR_\infty$.
As we discussed before, since $a_m^B$ has order $m$ sum rules, we have $\sd_{a_m^B} (\pp_{a_m^B}^{(j)})=2^{m-1-j}\pp_{a_m^B}^{(j)}$ for all $j=0,\ldots,m-1$, where $\pp_{a^B_m}$ is defined in \eqref{pa} with $a,m$ being replaced by $a^B_m$ and $m-1$, respectively. Define vector polynomials by
\[
\vec{\pp}_j:=(\pp_{a_m^B}^{(m-1-j)},0),
\quad j=0,\ldots,m-1 \quad \mbox{and}\quad
\vec{\pp}_{m+k}:=(0, \pp_{a_n^B}^{(n-1-k)}),\qquad k=0,\ldots,n-1.
\]
Then we have $\sd_a \vec{\pp}_j=2^{-j} \vec{\pp}_j$ for all $j=0,\ldots,m+n-1$.
Note that $\deg(\vec{\pp}_{m+k})=k\ne m+k$ for all $k=0,\ldots,n-1$.
Using \cref{cor:sr} and noting that $\deg(\vec{\pp}_j)=j$ for all $j=0,\ldots,m-1$, we see that $a$ has order $m$ sum rules. However, the mask $a$ cannot have order $m+n$ sum rules.
Suppose not. By \cref{cor:sr}, there exists $\vec{\pp}\in (\PL_{m+n-1})^{1\times 2}$ such that $\deg(\vec{\pp})=m+n-1$ and $\sd_a \vec{\pp}=2^{1-m-n} \vec{\pp}$.
Write $\vec{\pp}=(\pp_1,\pp_2)$. Since $\deg(\vec{\pp})=m+n-1$, without loss of generality, we assume  $\deg(\pp_1)=m+n-1$.
Then $\sd_a \vec{\pp}=2^{1-m-n} \vec{\pp}$ must imply $\sd_{a^B_m} \pp_1=2^{1-m-n} \pp_1$. Since $\deg(\pp_1)=m+n-1$, by \cref{cor:sr}, the mask $a_m^B$ must have order $m+n$ sum rules, which is a contradiction to $\wh{a_m^B}(\xi)=2^{-m} (1+e^{-i\xi})^m$.
Hence, the mask $a$ cannot have order $m+n-1$ sum rules, even though its vector subdivision operator $\sd_a$ has eigenvalues $2^{-j}, j=0,\ldots,m+n-1$ with all eigenvectors being vector polynomials in $\PR_\infty$.

Now suppose that we remove the factor $2^{-m}$ before $\wh{a^B_n}$ in the definition of the matrix-valued mask $a$. Define $\vgu_1,\vgu_2\in \lrs{0}{1}{2}$ by
\[
\wh{\vgu_1}(\xi):=(1/\wh{\varphi_{a^B_m}}(\xi),0)+\bo(|\xi|^m)
\quad \mbox{and}\quad
\wh{\vgu_2}(\xi):=(0,1/\wh{\varphi_{a^B_n}}(\xi))+\bo(|\xi|^n),\quad \xi\to 0,
\]
where $\wh{\varphi_{a^B_m}}(\xi):=\prod_{j=1}^\infty \wh{a^B_m}(2^{-j}\xi)$. Note that $\wh{\vgu_1}(0)=(1,0)\ne 0$ and $\wh{\vgu_2}(0)=(0,1)\ne 0$.
Then we can easily check that the mask $a$ (after dropping the factor $2^{-m}$) has order $m$ sum rules with respect to $\vgu_1\in \lrs{0}{1}{2}$, but the same mask $a$ also has order $n$ sum rules with respect to $\vgu_2\in \lrs{0}{1}{2}$.
}\end{example}

\section{Characterization of Matrix-valued Masks for Hermite Subdivision Schemes}
\label{sec:hsd}

To characterize masks for convergent Hermite subdivision schemes,
in this section we shall first study some necessary conditions for convergent Hermite subdivision schemes of order $r$ by linking them to vector cascade algorithms and refinable vector functions in wavelet theory.
This allows us to characterize the matrix-valued masks for convergent Hermite subdivision scheme in \cref{thm:hsd:mask} for which we shall provide a proof in this section. The characterization of convergence of a Hermite subdivision scheme stated in \cref{thm:hsd:converg} will be proved later in \cref{sec:converg}.

To link Hermite subdivision schemes with vector cascade algorithms and refinable vector functions, we recall that
the \emph{refinement operator}
$\cd_a: (\Lp{p})^r \rightarrow (\Lp{p})^r$
is defined in \eqref{cd}, that is,
$\cd_a f:=2\sum_{k\in \Z} a(k) f(2\cdot-k)$ for $f\in (\Lp{p})^r$.
By induction, we can easily observe that
\be \label{cdn}
\cd_a^n f=\sum_{k\in \Z} (\sd_a^n (\td I_r))(k) f(2^n\cdot-k)\quad \mbox{or equivalently},\quad \wh{\cd_a^n f}(\xi)=\wh{a_n}(2^{-n}\xi)\wh{f}(2^{-n}\xi),
\ee
where $\wh{a_n}(\xi)=\wh{a}(2^{n-1}\xi)\cdots \wh{a}(2\xi)\wh{a}(\xi)$ as defined in \eqref{an}.
The cascade algorithm produces a sequence $\{\cd_a^n f\}_{n=1}^\infty$ of iteratively generated vector functions, which may converge to a limiting vector function in some function spaces.

Let us first prove a simple fact, which is critical for us to link Hermite subdivision schemes with cascade algorithms and refinable vector functions. Because we are not aware of any explicit proof of this simple fact for Hermite subdivision schemes, we provide a proof here.

\begin{prop}\label{prop:hsd:phi}
Let $r\in \N$ and $m\in \NN$ with $m\ge r-1$.
Assume that the Hermite subdivision scheme of order $r$ associated with a mask $a\in \lrs{0}{r}{r}$ is convergent with limiting functions in $\CH{m}$ (see \cref{def:hsd}). Let the compactly supported vector function $\phi\in (\CH{m})^r$ be defined as the limiting vector function through \eqref{sdn:phi} using the initial data $w_0=\td I_r\in \lrs{0}{r}{r}$.
Then the compactly supported vector function $\phi$ must satisfy the following refinement equation:
\be \label{refeq}
\phi=2\sum_{k\in \Z} a(k) \phi(2\cdot-k)\quad \mbox{or equivalently},\quad
\wh{\phi}(2\xi)=\wh{a}(\xi)\wh{\phi}(\xi).
\ee
Moreover, $\wh{\phi}(\xi)=\lim_{n\to \infty} \big[\prod_{j=1}^n \wh{a}(2^{-j}\xi)\big ]e_1$ for $\xi\in \R$.
\end{prop}

\bp Define $h(x):=\max(1-|x|,0)$ to be the centered hat function.
For $n\in \N$, we define a sequence of vector functions $f_n, n\in \N$ by
\[
f_n:=\sum_{k\in \Z} (\sd_a^n (\td I_r))(k)f_0(2^{n}\cdot-k)=
2^n \sum_{k\in \Z} a_n(k) f_0(2^n\cdot-k),
\]
where $f_0:=(h,0,\ldots,0)^\tp$ and $a_n$ is defined in \eqref{an}. That is,
$f_n =\cd_a^n f_0$. We now prove that $\lim_{n\to \infty}\|f_n-\phi\|_{(\CH{})^r}=0$. Since $h(k)=\td(k)$ for all $k\in \Z$,
we notice that
$f_n(2^{-n}k)=2^n a_n(k)e_1$ for all $k\in \Z$ and $n\in \N$. Therefore,
by $f_0=e_1 h$, we obtain
\[
f_n=\sum_{k\in \Z} (\sd_a^n (\td I_r))(k)e_1 h(2^n\cdot-k)=
2^n \sum_{k\in \Z} a_n(k) e_1 h(2^n\cdot-k).
\]
Because \eqref{sdn:phi} holds and $\sD^{-n}e_1=e_1$, in particular, \eqref{sdn:phi} must hold for its first column, that is,
\be \label{sdn:phi:e1}
\lim_{n\to \infty} \|2^n a_n(\cdot)e_1-\phi(2^{-n}\cdot)\|_{\lrs{\infty}{r}{r}}=0,
\ee
where we used
$\sd_a^n (\td I_r)=2^n a_n$ by \eqref{an}.
Define a sequence of vector functions $\psi_n, n\in \N$ by
\[
\psi_n:=\sum_{k\in \Z} \phi(2^{-n}k)h(2^n\cdot-k).
\]
Then $\psi_n(2^{-n}k)=\phi(2^{-n}k)$ for all $k\in \Z$ and $n\in \N$. Since $\phi$ is a compactly supported continuous vector function, $\phi$ must be uniformly continuous and consequently,
\be \label{psin}
\lim_{n\to \infty} \|\psi_n-\phi\|_{(\CH{})^r}=0.
\ee
Since $\sum_{k\in\Z} |h(x+k)|\le 1$ for all $x\in \R$, we conclude that
\begin{align*}
\|f_n-\psi_n\|_{(\CH{})^r}
&=\left\| \sum_{k\in \Z} \left( 2^n a_n(k)e_1-\phi(2^{-n}k)\right) h(2^n\cdot-k) \right\|_{(\CH{})^r}\\
&\le \|2^n a_n(\cdot)e_1-\phi(2^{-n}\cdot)\|_{\lrs{\infty}{r}{r}}.
\end{align*}
As a consequence, we have
\[
\|f_n-\phi\|_{(\CH{})^r}\le
\|f_n-\psi_n\|_{(\CH{})^r}
+\|\psi_n-\phi\|_{(\CH{})^r}
\le \|2^n a_n(\cdot)e_1-\phi(2^{-n}\cdot)\|_{\lrs{\infty}{r}{r}}+
\|\psi_n-\phi\|_{(\CH{})^r}.
\]
Now it follows directly from \eqref{sdn:phi:e1} and \eqref{psin} that $\lim_{n\to \infty} \|f_n-\phi\|_{(\CH{})^r}=0$. That is,
$\lim_{n\to \infty} \|\cd_a^n f_0-\phi\|_{(\CH{})^r}=0$.
On the other hand, by induction on $n$ we observe that
\[
\wh{f_n}(\xi)=\wh{a_n}(2^{-n}\xi)\wh{f_0}(2^{-n}\xi)
=\wh{a}(2^{-1}\xi)\wh{a_{n-1}}(2^{-n}\xi)\wh{f_0}(2^{-n}\xi)
=\wh{a}(\xi/2)\wh{f_{n-1}}(\xi/2).
\]
That is, we proved $\wh{f_n}(2\xi)=\wh{a}(\xi)\wh{f_{n-1}}(\xi)$, which is equivalent to
\[
f_n=\cd_a f_{n-1}=2\sum_{k\in\Z} a(k) f_{n-1}(2\cdot-k).
\]
Since $a\in \lrs{0}{r}{r}$ has finite support, we can easily deduce that all the vector functions $f_n$ and $\phi$ are supported inside $\fs(a)$. Because $\lim_{n\to \infty} \|f_n-\phi\|_{(\CH{})^r}=0$, we now conclude from the above identity that \eqref{refeq} must hold for $\phi$.
Moreover, we have $\wh{\phi}(\xi)=\lim_{n\to \infty} \wh{f_n}(\xi)=\lim_{n\to \infty} \wh{a_n}(2^{-n}\xi) e_1\wh{h}(2^{-n}\xi)$. Since $\lim_{n\to\infty} \wh{h}(2^{-n}\xi)=\wh{h}(0)=1$, we conclude that
\[
\wh{\phi}(\xi)=\lim_{n\to \infty} \wh{a_n}(2^{-n}\xi)e_1=\lim_{n\to \infty} \Big[\prod_{j=1}^n \wh{a}(2^{-j}\xi)\Big ]e_1
\]
for $\xi\in \R$. This completes the proof.
\ep

We shall use the existence of a refinable vector function in \cref{prop:hsd:phi} to link Hermite subdivision schemes with cascade algorithms and refinable vector functions. To do so, we need the following special case of a technical result from
\cite[Proposition~5.6.2]{hanbook} and \cite[Section~3]{han03}.

\begin{prop}\label{prop:phi}
(\cite[Proposition~5.6.2]{hanbook} and \cite[Section~3]{han03})
Let $a\in \lrs{0}{r}{r}$ and $m\in \NN$.
Let $\phi=(\phi_1,\ldots,\phi_r)^\tp$ be a compactly supported vector function in $(\CH{m})^r$ satisfying the refinement equation in \eqref{refeq} and the condition in \eqref{cond:phi}.
%
%
Then the following statements hold:
\begin{enumerate}
\item[(1)]
$1$ is a simple eigenvalue of $\wh{a}(0)$ and all other eigenvalues of $\wh{a}(0)$ are less than $2^{-m}$ in modulus.

\item[(2)] There exists $\vgu_a\in \lrs{0}{1}{r}$ such that $\wh{\vgu_a}(0)\wh{\phi}(0)=1$ and the mask $a$ has order $m+1$ sum rules with respect to the matching filter $\vgu_a$, i.e., \eqref{sr} holds.

\item[(3)] All the vectors $\wh{\vgu_a}^{(j)}(0), j=0,\ldots,m$ are uniquely determined through the recursive formula: $\wh{\vgu_a}(0)\wh{a}(0)=\wh{\vgu_a}(0)$ with the normalization condition $\wh{\vgu_a}(0)\wh{\phi}(0)=1$ and
\be \label{vgua:value}
\wh{\vgu_a}^{(j)}(0)=[1-2^j\wh{a}(0)]^{-1}\sum_{k=0}^{j-1} \frac{2^k j!}{k!(j-k)!} \wh{\vgu_a}^{(k)}(0)\wh{a}^{(j-k)}(0),\qquad j=1,\ldots,m.
\ee

\item[(4)] Any polynomial $\pp\in \PL_m$ can be reproduced through $\pp=\sum_{k\in \Z}(\pp*\vgu_a)(k)\phi(\cdot-k)$ and
\be \label{phi:poly}
\wh{\vgu_a}(\xi)\wh{\phi}(\xi+2\pi k)=\td(k)+\bo(|\xi|^{m+1}),\qquad \xi\to 0, \; \forall\, k\in \Z.
\ee
\end{enumerate}
\end{prop}

\bp
Since $\phi\in (\CH{m})^r$ and $\cd_a \phi=\phi$ by \eqref{refeq}, it is
trivial that $\lim_{n\to \infty} \|\cd_a^n \phi-\phi\|_{(\CH{m})^r}=0$.
Note that the condition in \eqref{cond:phi} prevents $\phi$ to be identically zero.
Now all the claims in items (1)--(4) follow directly from \cite[Proposition~5.6.2]{hanbook} by taking $f=\phi$. Also, see \cite[Section~3]{han03} for closely related results for multivariate refinable vector functions.
\ep

Using the above two results, we are now ready to prove \cref{thm:hsd:mask} characterizing all the masks for convergent Hermite subdivision schemes.

\bp[Proof of \cref{thm:hsd:mask}]
By \cref{prop:hsd:phi}, we know that $\phi\in (\CH{m})^r$ is a refinable vector function satisfying \eqref{refeq} and $\wh{\phi}(\xi)=\lim_{n\to \infty} \big[\prod_{j=1}^n \wh{a}(2^{-j}\xi)\big ]e_1$ for $\xi\in \R$.
By \cref{prop:phi}, all the claims in items (1)--(4) of \cref{prop:phi} hold.
By $\wh{\vgu_a}(0)\wh{\phi}(0)=1$ and
item (2) of \cref{prop:phi}, we have
\[
1=\wh{\vgu_a}(0)\wh{\phi}(0)
=\lim_{n\to \infty} \wh{\vgu_a}(0)[\wh{a}(0)]^n e_1=\lim_{n\to \infty} \wh{\vgu_a}(0) e_1= \wh{\vgu_a}(0) e_1.
\]
That is, we proved $\wh{\vgu_a}(0)e_1=1$.
By item (2) of \cref{prop:phi}, we deduce from item (5) of \cref{thm:sd:sr} that $\sd_a \vec{\pp}_m=2^{-m}\vec{\pp}_m$ with $\vec{\pp}_m:=\frac{(\cdot)^m}{m!}*\vgu_a$ and $\deg(\vec{\pp}_m)=m$.
Therefore, $\sd_a^n \vec{\pp}_m=2^{-mn}\vec{\pp}_m$ for all $n\in \N$. Take $w_0=\vec{\pp}_m$ as the initial vector sequence. Then the refinement data $\{w_n\}_{n=1}^\infty$ in \eqref{hsd:wn} must satisfy $w_n=\sd_a^n (\vec{\pp}_m) \sD^{-n}=2^{-mn}\vec{\pp}_m \sD^{-n}$.
Since the Hermite subdivision scheme of order $r$ associated with mask $a$ is convergent with limiting functions in $\CH{m}$, there exists a function $\eta\in \CH{m}$ such that \eqref{hsd:converg:2} holds, that is, for any $K>0$,
\be \label{hsd:converg:3}
\lim_{n\to \infty}
\max_{k\in [-2^nK,2^nK]\cap \Z} | 2^{-mn} \vec{\pp}_m(k) 2^{\ell n} e_{\ell+1} -\eta^{(\ell)}(2^{-n}k)|=0,
\qquad \forall\, \ell=0,\ldots,r-1,
\ee
where we used
$\sD^{-n} e_{\ell+1}=2^{\ell n} e_{\ell+1}$.
On the other hand, by \cref{lem:conv} and $\vec{\pp}_m=\frac{(\cdot)^m}{m!}*\vgu_a$, we have
\[
2^{(\ell-m)n} \vec{\pp}_m(\cdot) e_{\ell+1}
=2^{(\ell-m)n}
\sum_{j=0}^m \frac{(-i)^j}{j!}
\frac{(\cdot)^{m-j}}{(m-j)!} \wh{\vgu_a}^{(j)}(0)e_{\ell+1}
=\sum_{j=0}^m \frac{(-i)^j}{j!} \frac{(2^{-n}\cdot)^{m-j}}{(m-j)!}
\wh{\vgu_a}^{(j)}(0) e_{\ell+1} 2^{(\ell-j)n}.
\]
For $n\in \N$, we define functions
\be \label{psi:ell}
\psi_{\ell+1,n}(x):=\sum_{j=0}^m \frac{(-i)^j}{j!} \frac{x^{m-j}}{(m-j)!}
\wh{\vgu_a}^{(j)}(0)e_{\ell+1} 2^{(\ell-j)n},\qquad \ell=0,\ldots,r-1.
\ee
Then
$2^{(\ell-m)n}\vec{\pp}_m(k) e_{\ell+1}=\psi_{\ell+1,n}(2^{-n}k)$ for all $k\in \Z$ and $n\in \N$. Therefore,
\eqref{hsd:converg:3} implies
\be \label{lim:psin}
\lim_{n\to\infty}
\psi_{\ell+1,n}(x)=\eta^{(\ell)}(x) \qquad \forall\, x\in \cup_{n=1}^\infty (2^{-n} \Z)\quad \mbox{and}\quad \ell=0,\ldots,r-1.
\ee
By the definition of $\psi_{\ell+1,n}$ in \eqref{psi:ell} and
\[
\lim_{n\to \infty} 2^{(\ell-j)n}=\begin{cases}
0, &\text{if $j>\ell$},\\
1, &\text{if $j=\ell$},\\
\infty, &\text{if $j<\ell$},
\end{cases}
\]
the existence of the limits in \eqref{lim:psin} forces
\be \label{hsd:vgua:eq1}
\wh{\vgu_a}^{(j)}(0)e_{\ell+1}=0,\qquad \forall\, 0\le j<\ell\le r-1
\ee
and
\be \label{hsd:vgua:eq2}
\eta^{(\ell)}(x)=\lim_{n\to \infty} \psi_{\ell+1,n}(x)=
\frac{(-i)^\ell}{\ell!} \frac{x^{m-\ell}}{(m-\ell)!} \wh{\vgu_a}^{(\ell)}(0) e_{\ell+1}=:\pq_{\ell+1}(x),\qquad \ell=0,\ldots,r-1.
\ee
Since $\cup_{n=1}^\infty (2^{-n}\Z)$ is dense in $\R$, the above identity in \eqref{lim:psin} must hold for all $x\in \R$ and therefore,
\[
\frac{(-i)^\ell}{\ell!} \frac{x^{m-\ell}}{(m-\ell)!} \wh{\vgu_a}^{(\ell)}(0) e_{\ell+1}=
\pq_{\ell+1}(x)=[\eta(x)]^{(\ell)}=\pq_1^{(\ell)}(x)
=\frac{x^{m-\ell}}{(m-\ell)!} \wh{\vgu_a}(0)e_1.
\]
Because we proved $\wh{\vgu_a}(0)e_1=1$, the above identity becomes
\be \label{hsd:vgua:eq3}
\wh{\vgu_a}^{(\ell)}(0) e_{\ell+1}=\ell!i^\ell,\qquad \forall\, \ell=0,\ldots,r-1.
\ee
Note that \eqref{hsd:vgua:eq1} and \eqref{hsd:vgua:eq3} together are equivalent to
\be \label{hsd:vgua:0}
\wh{\vgu_a}(\xi)e_{\ell+1}=(i\xi)^\ell +\bo(|\xi|^{\ell+1}),\qquad \xi \to 0, \; \forall\; \ell=0,\ldots,r-1,
\ee
which is just \eqref{vgu:hermite:0}.
In particular,
$\wh{\vgu_a}(0)=(1,0,\ldots,0)$ and thus $e_1^\tp \wh{\phi}(0)=\wh{\vgu_a}(0)\wh{\phi}(0)=1$.
Now \eqref{hsd:vgua:0} is equivalent to \eqref{vgu:hermite} for some $c_1,\ldots,c_r\in \lp{0}$ with $\wh{c_1}(0)=\cdots=\wh{c_r}(0)=1$.
This proves all the claims in \cref{thm:hsd:mask}.
\ep

Motivated by \cref{thm:hsd:mask}, for Hermite subdivision schemes, we define

\begin{definition}\label{def:hmask}
{\rm Let $r\in \N$ and $m\in \NN$ with $m\ge r-1$.
We say that $a\in \lrs{0}{r}{r}$ is \emph{a Hermite mask of accuracy order $m+1$} if the mask $a$ has order $m+1$ sum rules with respect to some $\vgu_a\in \lrs{0}{1}{r}$ such that
\eqref{vgu:hermite:0} or
\eqref{vgu:hermite} is satisfied for some $c_1,\ldots,c_r\in \lp{0}$ with $\wh{c_1}(0)=\cdots=\wh{c_r}(0)=1$.
}
\end{definition}

By \cref{prop:phi} and Theorem~\ref{thm:hsd:mask},
for a convergent Hermite subdivision scheme associated with a mask $a\in \lrs{0}{r}{r}$ with limiting functions in $\CH{m}$, items (1) and (3) of \cref{prop:phi} must hold. Hence,
every finitely supported Hermite mask $a\in \lrs{0}{r}{r}$ of accuracy order $m+1$ can be easily obtained by solving the following linear system:
$\wh{\vgu_a}(0)=(1,0,\ldots,0)$ and
\be \label{hmask}
\wh{\vgu_a}^{(j)}(0)=\sum_{k=0}^{j} \frac{2^k j!}{k!(j-k)!} \wh{\vgu_a}^{(k)}(0)\wh{a}^{(j-k)}(0),\qquad j=1,\ldots,m,
\ee
where the vectors $\wh{\vgu_a}^{(j)}(0), j=1,\ldots,m$ are given in \eqref{vgu:hermite}, for finding the unknown coefficients of a mask $a\in \lrs{0}{r}{r}$ and $\wh{c_\ell}^{(j)}(0), 1\le \ell \le r$ and $1\le j\le m+1-\ell$ with $\wh{c_1}(0)=\cdots=\wh{c_r}(0)=1$.

For an interpolatory Hermite subdivision scheme of order $r$ associated with a mask $a\in \lrs{0}{r}{r}$, its mask $a$ must satisfy the interpolation condition in \eqref{sd:int}. If in addition $a$ is a Hermite mask of accuracy order $m+1$, then it is easy to deduce from \eqref{sd:int} and \eqref{sr:2} that the sequence $\vgu_a$ satisfies \eqref{vgu:hermite}
if and only if $\wh{c_\ell}(\xi)=1+\bo(|\xi|^{m+1})$ as $\xi \to 0$ for all $\ell=1,\ldots,r$, that is, $\vgu_a$ is just $\vgu_H\in \lrs{0}{1}{r}$ satisfying
\be \label{vguH}
\wh{\vgu_H}(\xi)=(1,i\xi,\ldots,(i\xi)^{r-1})
+\bo(|\xi|^{m+1}),\qquad \xi\to 0,
\ee
which agrees with \cite[Lemma~4.1]{han01}, \cite[Proposition~5.3]{han03} and \cite[Lemma~6.2.5]{hanbook} for interpolatory Hermite subdivision schemes. Moreover, by \cref{thm:sd:sr}, we have $\sd_a \vec{\pp}_j=2^{-j}\vec{\pp}_j$ with $\vec{\pp}_j:=\pp_j*\vgu_H=(\pp_j,\pp_j',\ldots,\pp_j^{(r-1)})$ and $\pp_j(x):=\frac{x^j}{j!}$ for $j=0,\ldots,m$. This eigenvalue condition is called the spectral condition in \cite{dm09} for studying Hermite subdivision schemes.

As we discussed before, for any initial input data $w_0\in \lrs{0}{1}{r}$, its Hermite refinements $\{w_n\}_{n=1}^\infty$ defined in \eqref{hsd:wn} converge to $w_0*\phi:=\sum_{k\in \Z} w_0(k)\phi(\cdot-k)$. It is natural to require that $w_0*\phi$ be identically zero only if $w_0$ is identically zero. This is equivalent to saying that the integer shifts of $\phi$ are \emph{linearly independent}, which is further equivalent to
\be \label{li}
\mspan\{\wh{\phi}(\xi+2\pi k) \setsp k\in \Z\}=\C^r,\qquad \forall\, \xi\in \C.
\ee
For a simple proof of the above equivalence, see \cite{han16} and references therein. Obviously, \eqref{li} implies the condition in \eqref{cond:phi} of \cref{thm:hsd:mask} and the stability condition in \eqref{stability} of Theorem~\ref{thm:hsd:converg}. On the other hand, as explained in \cite[Theorem~2.3]{han01} and \cite[Theorem~5.5.4]{hanbook}, the condition in \eqref{cond:phi} is often necessary to guarantee the sum rules condition in item (2) of \cref{prop:phi}. Hence, \eqref{cond:phi} is probably the weakest condition to avoid degenerate Hermite subdivision schemes.

\section{Factorization and Convergence of Hermite Subdivision Schemes}
\label{sec:converg}

In this section we shall introduce the normal form and factorization of a general matrix-valued mask and then we shall characterize convergence of Hermite subdivision schemes by proving \cref{thm:hsd:converg}.

The analysis of convergence of general vector subdivision schemes and smoothness of refinable vector functions has been widely known to be much more difficult and complicated than their scalar counterparts. This is largely because general vector subdivision schemes and refinable vector functions employ matrix-valued masks and as we already discussed in \cref{sec:sd}, the notion of sum rules for matrix-valued masks is much more involved than their scalar counterparts.

To facilitate the study of refinable vector functions,
the normal form of a matrix-valued mask has been initially introduced in \cite[Theorems~2.2 and~2.3]{hm03} for studying one-dimensional dual multiframelets and in \cite[Proposition~2.4]{han03} for studying convergence of multivariate vector cascade algorithms and smoothness of refinable vector functions. The normal form of matrix-valued masks has been further developed in \cite[Theorem~2.1]{han09} and \cite[Theorem~5.1]{han10mc} for studying balancing properties of multiwavelets and multiframelets, which are derived from refinable vector functions.

In order to introduce the normal form of a matrix-valued mask and study the convergence of Hermite subdivision schemes, let us first introduce some definitions.
For $\vgu\in \lrs{0}{1}{r}$ and $m\in \NN$, recall that
$\PR_{m,\vgu}:=\{\pp*\vgu \setsp \pp\in\PL_m\}$. For analyzing convergence of vector cascade algorithms and refinable vector functions, a closely related dual space $\PV_{m,v}$ is defined to be
\be \label{Vvm}
\PV_{m,\vgu}:=\{u\in (\lp{0})^r \setsp \wh{\vgu}(\xi)\wh{u}(\xi)=\bo(|\xi|^{m+1}),\quad \xi\to 0\}.
\ee
Note that $\PV_{m,\vgu}$ is shift-invariant, that is, $u\in \PV_{m,\vgu}$ implies $u(\cdot-k)\in \PV_{m,\vgu}$ for all $k\in \Z$.
It is important to notice that
we removed both complex conjugate and transpose in the definition of
the space $\PV_{m,\vgu}$ in \cite{han03,hanbook}, which is a subspace of $\lrs{0}{1}{r}$ (instead of $(\lp{0})^r$ here) and is just $\{ \ol{u(-\cdot)}^\tp \setsp u\in \PV_{m,\vgu}\}$ with $\PV_{m,\vgu}$ here defined in \eqref{Vvm}.
The definition in \eqref{Vvm} makes the presentation simple.

For $U\in \lrs{0}{r}{r}$, we say that $U$ (or $\wh{U}$) is \emph{strongly invertible} if $\det(\wh{U})$ is a nonzero monomial, in other words, $(\wh{U}(\xi))^{-1}$ is a matrix of $2\pi$-periodic trigonometric polynomials.
If $U\in \lrs{0}{r}{r}$ is strongly invertible, then we can define $U^{-1}\in \lrs{0}{r}{r}$ by $\wh{U^{-1}}(\xi):=(\wh{U}(\xi))^{-1}$. Note that $U*U^{-1}=U^{-1}*U=\td I_r$.
For $u\in \lrs{0}{s}{r}$ and $n\in \NN$, we define $\nabla u:=u-u(\cdot-1)$ and $\nabla^n u:=\nabla^{n-1} (\nabla u)$, that is, $\wh{\nabla^n u}(\xi):=(1-e^{-i\xi})^n\wh{u}(\xi)$.

Mainly following \cite[Theorem~5.6.4]{hanbook}, we now state the result on the normal form of a matrix-valued mask, for which we give a sketch of proof here.

\begin{theorem}\label{thm:nf}
Let $m\in \NN$ and $a\in \lrs{0}{r}{r}$ such that there is $\vgu_a\in \lrs{0}{1}{r}$ satisfying
\be \label{cond:vgua}
\wh{\vgu_a}(0)\ne 0 \quad \mbox{and}\quad
\wh{\vgu_a}(2\xi)\wh{a}(\xi)=\wh{\vgu_a}(\xi)+\bo(|\xi|^{m+1}),\qquad \xi \to 0.
\ee
Then there exists a strongly invertible sequence $U\in \lrs{0}{r}{r}$ such that
\be \label{nf:vgua}
\wh{\mathring{\vgu}}(\xi):=\wh{\vgu_a}(\xi)\wh{U}(\xi)
=(1+\bo(|\xi|),\bo(|\xi|^{m+1})),\ldots, \bo(|\xi|^{m+1}), \qquad \xi \to 0
\ee
and the following statements hold:
\begin{enumerate}
\item[(1)] If $\phi$ is an $r\times 1$ vector of compactly supported distributions satisfying $\wh{\phi}(2\xi)=\wh{a}(\xi)\wh{\phi}(\xi)$, define
\be \label{nf:mask:phi}
\wh{\mathring{a}}(\xi):=(\wh{U}(2\xi))^{-1} \wh{a}(\xi) \wh{U}(\xi) \quad \mbox{and}\quad
\wh{\mathring{\phi}}(\xi):= (\wh{U}(\xi))^{-1}\wh{\phi}(\xi),
\ee
then $\mathring{\phi}$ is an $r\times 1$ vector of compactly supported distributions satisfying $\wh{\mathring{\phi}}(2\xi)=\wh{\mathring{a}}(\xi)\wh{\mathring{\phi}}(\xi)$ and its associated mask $\mathring{a}$ must be finitely supported, i.e., $\mathring{a}\in \lrs{0}{r}{r}$.
\item[(2)] The finitely supported mask $\mathring{a}$ must take the following normal (or canonical) form:
\be \label{nf:mask}
\left[
\begin{matrix}
a_{1,1} &a_{1,2}\\
a_{2,1} &a_{2,2}\end{matrix}\right]
\quad \mbox{with}\quad \wh{a_{1,1}}(0)=1, \quad \wh{a_{1,2}}(\xi)=\bo(|\xi|^{m+1}),\quad \xi \to 0,
\ee
where $a_{1,1}\in \lp{0}$, $a_{1,2}\in \lrs{0}{1}{(r-1)}$, $a_{2,1}\in \lrs{0}{(r-1)}{1}$, and $a_{2,2}\in \lrs{0}{(r-1)}{(r-1)}$.

\item[(3)] The mask $a$ has order $m+1$ sum rules with respect to $\vgu_a$ if and only if $\mathring{a}$ has order $m+1$ sum rules with respect to $\mathring{\vgu}$, which is further equivalent to that $\mathring{a}$ satisfies \eqref{nf:mask} and
\be \label{nf:mask:sr}
\wh{a_{1,1}}(\xi+\pi)=\bo(|\xi|^{m+1}), \qquad \wh{a_{1,2}}(\xi+\pi)=\bo(|\xi|^{m+1}),\qquad \xi \to 0.
\ee
\item[(4)] $\PR_{m,\vgu_a}=\PR_{m,\mathring{\vgu}}*U^{-1}:=
    \{u*U^{-1} \setsp u\in \PR_{m,\mathring{\vgu}}\}$ and
 $\PV_{m,\vgu_a}=U*\PV_{m,\mathring{\vgu}}$.
\item[(5)] $\PV_{m,\vgu_a}$ is generated by $\PB_{m,\vgu_a}$, i.e., $\PV_{m,\vgu_a}=\mspan\{u(\cdot-k) \setsp u\in \PB_{m,\vgu_a},k\in \Z\}$, where
\be \label{BV}
\PB_{m,\vgu_a}:=U*\PB_{m,\mathring{\vgu}} \quad \mbox{with}\quad
\PB_{m,\mathring{\vgu}}:=\{ \nabla^{m+1} \td e_1, \td e_2,\ldots, \td e_r\}.
\ee

\item[(6)] For the case $r>1$, we can further choose such a strongly invertible sequence $U\in \lrs{0}{r}{r}$ satisfying all the above claims in items (1)--(5) with the additional properties
\be \label{nf:extra}
\wh{\mathring{\vgu}}(\xi):=\wh{\vgu_a}(\xi)\wh{U}(\xi)=(1,0,\ldots,0)+\bo(|\xi|^{m+1})
\quad \mbox{and}\quad
\wh{a}(\xi)=1+\bo(|\xi|^{m+1}),\qquad \xi\to 0.
\ee
\end{enumerate}
\end{theorem}

\bp We sketch the main idea of proof. The reader is referred to \cite{han03,han09,han10mc,hanbook,hm03} for details.
The key idea for proving all the claims in \cref{thm:nf} is very simple: We can always construct a strongly invertible sequence $U\in \lrs{0}{r}{r}$ such that \eqref{nf:vgua} holds.
To do so, we write $(\vgu_1,\ldots,\vgu_r):=\vgu_a$ by listing the entries of the vector $\vgu_a$.
Because $\wh{\vgu_a}(0)\ne 0$, without loss of generality we can assume $\wh{\vgu_1}(0)\ne 0$, otherwise we perform a permutation on the entries of $\vgu_a$. Since $\wh{\vgu_1}(0)\ne 0$, we can easily find $u_2,\ldots,u_r\in \lp{0}$ such that
\[
\wh{u_\ell}(\xi)=\wh{\vgu_\ell}(\xi)/\wh{\vgu_1}(\xi)+\bo(|\xi|^{m+1}),\qquad \xi \to 0, \ell=2,\ldots,r.
\]
Now we can define a matrix-valued sequence $U\in \lrs{0}{r}{r}$ by
\be \label{U}
U=\frac{1}{\wh{\vgu_1}(0)}\left[ \begin{matrix}
\td &-u_2 &\cdots &-u_r\\
0 &\td &\cdots &0\\
\vdots &\vdots &\ddots &\vdots\\
0 &0 &0 &\td\end{matrix}\right], \quad
\mbox{i.e.},
\quad
\wh{U}(\xi):=\frac{1}{\wh{\vgu_1}(0)}\left[ \begin{matrix}
1 &-\wh{u_2} &\cdots &-\wh{u_r}(\xi)\\
0 &1 &\cdots &0\\
\vdots &\vdots &\ddots &\vdots\\
0 &0 &0 &1\end{matrix}\right].
\ee
Since $\det(\wh{U}(\xi))=1/\wh{\vgu_1}(0)\ne 0$, the sequence $U$ is clearly strongly invertible and
\eqref{nf:vgua} holds.
We can gain the extra property in item (6) for $r>1$ by using a technique on matrices. See \cite[Theorem~2.1]{han09}
for more details for proving item (6).

Now all the claims in \cref{thm:nf} can be easily and directly verified.
\ep

Because $U$ is strongly invertible, the mask $a$ and its refinable vector function $\phi$ in \cref{thm:nf} can be equivalently transformed through \eqref{nf:mask:phi} into the new mask
$\mathring{a}$ and a new refinable vector function $\mathring{\phi}$. Due to the normal form of the new mask in \eqref{nf:mask} and \eqref{nf:mask:sr}, almost all analysis techniques for scalar masks and scalar refinable functions can be applied to matrix-valued masks and refinable vector functions. The normal form of matrix-valued masks in \eqref{nf:mask} greatly facilitates the study of vector subdivision schemes, refinable vector functions and multiwavelets. See \cite{han03,han09,han10mc,hanbook,hm03} for more details.
Moreover, the normal form in \cref{thm:nf} and \cref{thm:hsd:mask,thm:hsd:converg} together can transform any matrix-valued mask into a Hermite mask.

In order to prove \cref{thm:hsd:converg}, we need to recall the definition of a technical quantity $\sm_\infty(a)$ from \cite{han03,hanbook}.
Let $m$ be the largest possible integer such that items (1) and (2) of \cref{prop:phi} are satisfied (i.e., we take the highest order $m+1$ of sum rules).
For $1\le p\le \infty$, we define
\be \label{sma}
\sm_p(a):=\frac{1}{p}-\log_2 \rho_{m+1}(a,\vgu_a)_p
\ee
where
\be \label{rho}
\rho_{m+1}(a,\vgu_a)_p:=2\max\{\limsup_{n\to \infty} \|a_n*u \|_{(\lp{p})^r}^{1/n} \setsp u\in \PB_{m,\vgu_a}\},
\ee
where $a_n:=2^{-n} \sd^n (\td I_r)$ in \eqref{an}, i.e., $\wh{a_n}(\xi):=\wh{a}(2^{n-1}\xi)\cdots\wh{a}(2\xi)\wh{a}(\xi)$, and
$\PB_{m,\vgu_a}\subseteq \PV_{m,\vgu_a}$ generates $\PV_{m,\vgu_a}$, i.e., $\mspan\{u(\cdot-k)\setsp u\in \PB_{m,\vgu_a}\}=\PV_{m,\vgu_a}$.
The normal form of a matrix-valued mask and the technical quantity $\sm_p(a)$ play critical roles in studying convergence of vector cascade algorithms, smoothness of refinable vector functions, multiwavelets and framelets, see \cite{han03,hanbook} and references therein for details.

To provide an example to
demonstrate the advantages of the normal form of matrix-valued masks, we employ it here to factorize matrix-valued masks. Suppose that a mask $a\in \lrs{0}{r}{r}$ has order $m+1$ sum rules with respect to $\vgu_a\in \lrs{0}{1}{r}$. Then the new mask $\mathring{a}$ in \cref{thm:nf} must satisfy \eqref{nf:mask} and \eqref{nf:mask:sr}. It follows easily from \eqref{nf:mask} and \eqref{nf:mask:sr} (e.g., see \cite[Theorem~5.8.3]{hanbook}) that  the new derived matrix-valued mask $b$ from the mask $\mathring{a}$, defined through
\[
\wh{b}(\xi):=
(\wh{D_{m+1}}(2\xi))^{-1}
\wh{\mathring{a}}(\xi) \wh{D_{m+1}}(\xi)
\quad \mbox{with}\quad
\wh{D_{m+1}}(\xi):=\mbox{diag}( (1-e^{-i\xi})^{m+1}, I_r),
\]
must be finitely supported. Consequently, we can factorize the matrix-valued mask $a$ as follows:
\be \label{a:V}
\wh{a}(\xi)=\wh{U}(2\xi) \wh{\mathring{a}}(\xi) (\wh{U}(\xi))^{-1}=
\wh{V}(2\xi) \wh{b}(\xi) (\wh{V}(\xi))^{-1}\quad \mbox{with}\quad
\wh{V}(\xi):=\wh{U}(\xi) \wh{D_{m+1}}(\xi).
\ee
Before preceding further, let us look at the matrix sequence $V$ first.
Write $U=[u_1,\ldots,u_r]$ with $u_1,\ldots,u_r$ being the column vectors of $U$. By the definition of $V$ in \eqref{a:V}, we can easily observe that $V=[\nabla^{m+1} u_1, u_2,\ldots, u_r]$. By the definition of the vector subdivision operator $\sd_a$ in \eqref{sd}, for $w\in \lrs{0}{s}{r}$, we have
\[
\wh{\sd_a w}(\xi)=2\wh{w}(2\xi) \wh{a}(\xi)=2 \wh{w}(2\xi) \wh{V}(2\xi) \wh{b}(\xi) (\wh{V}(\xi))^{-1}=\wh{\sd_{b} (w*V)}(\xi) (\wh{V}(\xi))^{-1}.
\]
That is, we arrive at the factorization of the vector subdivision operator for a matrix-valued mask $a$ with order $m+1$ sum rules as follows:
\be \label{factorization}
(\sd_a w)*V=\sd_b (w*V)\quad \mbox{with}\quad
V=[\nabla^{m+1}u_1, u_2,\ldots, u_r].
\ee
Consequently, we have
\[
(\sd_a^n w)*V=(\sd_a (\sd_a^{n-1} w))*V=
\sd_b ((\sd_a^{n-1} w)*V)=
\sd_b^n (w*V).
\]
Taking $w=\td I_r$ and noting that $\sd_a^n (\td I_r)=2^n a_n$, we obtain from the above identity that
\be \label{an:bn}
a_n*V=2^{-n} \sd_b^n V= (V(2^n\cdot))*b_n\quad \mbox{with}\quad
\wh{b_n}(\xi):=\wh{b}(2^{n-1}\xi)\cdots\wh{b}(2\xi)\wh{b}(\xi).
\ee
On the other hand, by the definition of $\PB_{m,\vgu_a}$ in \eqref{BV}, we notice that $U*(\nabla^{m+1}\td e_1)=\nabla^{m+1} u_1$ and $U*(\td e_\ell)=u_\ell$ for all $\ell=2,\ldots,r$. Hence, the columns of $V$ forms the basis elements in $\PB_{m,\vgu_a}$, which generates $\PV_{m,\vgu_a}$. In other words, by \eqref{an:bn} we reach
\[
\rho_{m+1}(a,\vgu_a)_p=2\max\{\limsup_{n\to \infty} \|a_n *V\|_{(\lp{p})^{r\times r}}^{1/n}\}=2\max\{\limsup_{n\to \infty} \|(V(2^n\cdot))*b_n\|_{(\lp{p})^{r\times r}}^{1/n}\}.
\]
Using a technical result in \cite{han03} or \cite[Theorem~5.8.3]{hanbook}, we conclude from the above identity that
\[
\rho_{m+1}(a,\vgu_a)_p=
\rho(b)_p:=2\max\{\limsup_{n\to \infty} \|b_n\|_{(\lp{p})^{r\times r}}^{1/n}\}.
\]
Therefore, the convergence of the vector subdivision scheme using $\sd_a$ can be characterized through its derived vector subdivision scheme using $\sd_b$ via the above identity. See \cite{han03,hanbook} for more details.

Let us demonstrate the normal form and factorization of Hermite masks for Hermite subdivision schemes.
Let $a\in \lrs{0}{r}{r}$ be a Hermite mask of accuracy order $m+1$ (see \cref{def:hmask} \eqref{vgu:hermite}). Hence, the mask $a$ has order $m+1$ sum rules with respect to $\vgu_a\in \lrs{0}{1}{r}$ satisfying \eqref{vgu:hermite} for some $c_1,\ldots,c_r\in \lp{0}$ with $\wh{c_1}(0)=\cdots=\wh{c_r}(0)=1$. That is, we have
\[
\wh{\vgu_a}(\xi)=\wh{c_1}(\xi)[1, i\xi \wh{c_2}(\xi)/\wh{c_1}(\xi),\ldots, (i\xi)^{r-1} \wh{c_r}(\xi)/\wh{c_1}(\xi)]+\bo(|\xi|^{m+1}),\quad \xi\to 0.
\]
Define $\wh{\eta}(\xi):=\frac{1-e^{-i\xi}}{i\xi}$, the Fourier transform of the characteristic function $\chi_{[0,1]}$. Then $\wh{\eta}$ is an infinitely differentiable function and $\wh{\eta}(0)=1$. Thus, for $\ell=2,\ldots,r$, since $\wh{c_1}(0)=\wh{\eta}(0)=1\ne 0$, we can take $d_{\ell}\in \lp{0}$ satisfying
\[
\wh{d_{\ell}}(\xi)=\frac{\wh{c_{\ell}}(\xi)}
{\wh{c_1}(\xi)(\wh{\eta}(\xi))^{\ell-1}}
+\bo(|\xi|^{m+1}),\qquad \xi\to 0.
\]
Consequently, for $\ell=2,\ldots,r$, we have
\[
\wh{\vgu_a}(\xi) e_{\ell}=
(i\xi)^{\ell-1} \frac{\wh{c_{\ell}}(\xi)}{\wh{c_1}(\xi)}+\bo(|\xi|^{m+1})
=(1-e^{-i\xi})^{\ell-1} \frac{\wh{c_{\ell}}(\xi)}{\wh{c_1}(\xi)
(\wh{\eta}(\xi))^{\ell-1}}+\bo(|\xi|^{m+1})
=(1-e^{-i\xi})^\ell \wh{d_{\ell}}(\xi)+\bo(|\xi|^{m+1})
\]
as $\xi \to 0$. That is, we have
\[
\wh{\vgu_a}(\xi)=\wh{c_1}(\xi)[1, (1-e^{i\xi})\wh{d_2}(\xi),\ldots, (1-e^{-i\xi})^{r-1} \wh{d_r}(\xi)]+\bo(|\xi|^{m+1}),\quad \xi \to 0.
\]
We can define a strongly invertible sequence $U\in \lrs{0}{r}{r}$ in \eqref{U} with $\wh{\vgu_1}(0)=\wh{c_1}(0)\ne 0$ and
\[
u_\ell:=\nabla^{\ell-1} d_\ell, \quad \mbox{that is}, \quad
\wh{u_\ell}(\xi):=(1-e^{-i\xi})^{\ell-1}\wh{d_\ell}(\xi),\qquad \ell=2,\ldots,r.
\]
More explicitly, the matrices $U$ in \eqref{U} and $V$ in \eqref{factorization}
can be explicitly written as
\[
U=\frac{1}{\wh{c_1}(0)}\left[ \begin{matrix}
 \td &-\nabla d_2 &\cdots &-\nabla^{r-1} d_r\\
0 &\td &\cdots &0\\
\vdots &\vdots &\ddots &\vdots\\
0 &0 &0 &\td\end{matrix}\right],
\qquad
V=\frac{1}{\wh{c_1}(0)}\left[ \begin{matrix}
 \nabla^{m+1}\td &-\nabla d_2 &\cdots &-\nabla^{r-1} d_r\\
0 &\td &\cdots &0\\
\vdots &\vdots &\ddots &\vdots\\
0 &0 &0 &\td\end{matrix}\right].
\]
Since $\det(\wh{U}(\xi))=1/\wh{c_1}(0)\ne 0$, the sequence $U$ is clearly strongly invertible and
\[
\wh{\mathring{\vgu}}(\xi):=\wh{\vgu_a}(\xi) \wh{U}(\xi)=[\wh{c_1}(\xi)/\wh{c_1}(0),0,\ldots,0]+\bo(|\xi|^{m+1}),\qquad \xi \to 0.
\]

To prove \cref{thm:hsd:converg}, we need the following auxiliary result.

\begin{lemma}\label{lem:to0}
Let $\psi$ be a compactly supported continuous function on $\R$. Let $u\in \lp{1}$ such that $\wh{u}(\xi)=1+(1-e^{-i\xi})\wh{c}(\xi)$ for some $c\in \lp{1}$, i.e., $u(k):=\td(k)+c(k)-c(k-1)$ for $k\in \Z$ (note that this condition is satisfied for any exponentially decaying sequence $u\in \lp{1}$ with $\wh{u}(0)=1$). Then $\lim_{n\to \infty} \|\psi(2^{-n}\cdot)*u-\psi(2^{-n}\cdot) \|_{\CH{}}=0$ and in particular, $\lim_{n\to \infty} \|\psi(2^{-n}\cdot)*u-\psi(2^{-n}\cdot)\|_{\lp{\infty}}=0$, where
\[
\psi(2^{-n}\cdot)*u:=\sum_{k\in \Z} \psi(2^{-n}(\cdot-k))u(k).
\]
\end{lemma}

\bp Observe that the Fourier transform of $\psi(2^{-n}\cdot)*u$ is $2^n \wh{\psi}(2^n\xi) \wh{u}(\xi)$ and
\[
2^n \wh{\psi}(2^n\xi) \wh{u}(\xi)=
2^n\wh{\psi}(2^n\xi)+
2^n \wh{\psi}(2^n\xi)(1-e^{-i\xi})\wh{c}(\xi)=
\wh{\psi(2^{-n}\cdot)}(\xi)+\wh{\psi_n*c}(\xi),
\]
where $\psi_n:=\psi(2^{-n}\cdot)-\psi(2^{-n}\cdot-2^{-n})$, that is, $\wh{\psi_n}(\xi):=2^n \wh{\psi}(2^n\xi) (1-e^{-i\xi})$. Consequently, we have $\psi(2^{-n}\cdot)*u=\psi(2^{-n}\cdot)+\psi_n*c$ and hence
\[
\|\psi(2^{-n}\cdot)*u-\psi(2^{-n}\cdot)\|_{\CH{}}
=\|\psi_n*c\|_{\CH{}}\le \|c\|_{\lp{1}} \|\psi_n\|_{\CH{}}
=\|c\|_{\lp{1}} \|\psi-\psi(\cdot-2^{-n})\|_{\CH{}}.
\]
Since $\psi$ is a compactly supported continuous function, $\psi$ must be uniformly continuous and thus $\lim_{n\to \infty} \|\psi-\psi(\cdot-2^{-n})\|_{\CH{}}=0$. This completes the proof.
\ep

We are now ready to prove \cref{thm:hsd:converg} on convergence of Hermite subdivision schemes.

\bp[Proof of \cref{thm:hsd:converg}]
(2)$\iff$(3) is established in \cite[Theorem~4.3]{han03}. Also, see \cite[Theorem~5.6.11]{hanbook} for more details and related results.

To prove (2)$\imply$(1), we need a suitable initial vector function $f$ satisfying \eqref{initialf}.
Because the existence of a compactly supported refinable Hermite interpolant $\theta\in (\CH{m})^r$ of order $r$ still remains open, we now modify a well-known Hermite interpolant $\theta:=(\theta_0,\ldots,\theta_m)^\tp \in (\CH{m})^{m+1}$ of order $m+1$ given below (e.g., see \cite[Proposition~6.2.2]{hanbook})
\[
\theta_{\ell}(x):=
\begin{cases}
(1-x)^{m+1} \frac{x^{\ell}}{\ell!}\sum_{j=0}^{m-\ell}
\frac{(m+j)!}{m!j!} x^j, &\text{$x\in [0,1]$},\\
(1+x)^{m+1} \frac{x^{\ell}}{\ell!}\sum_{j=0}^{m-\ell}
\frac{(m+j)!}{m!j!} (-x)^j, &\text{$x\in [-1,0)$},\\
0, &\text{$x\in \R\bs [-1,1]$},
\end{cases}
\]
for $\ell=0,\ldots,m$. Then $\theta\in (\CH{m})^{m+1}$ is a compactly supported refinable Hermite interpolant of order $m+1$ possessing the Hermite interpolation property:
\be \label{theta:hint}
\theta_{\ell}^{(j)}(k)=\td(\ell-j)\td(k),\qquad \forall\; \ell,j=0,\ldots,m\quad \mbox{and}\quad k\in \Z,
\ee
and by \cite[Corollary~5.2]{han03} or \cite[Theorem~6.2.3]{hanbook}, we have
\be \label{theta:approx}
(1,i\xi,\ldots, (i\xi)^m) \wh{\theta}(\xi+2\pi k)=\td(k)+\bo(|\xi|^{m+1}),\qquad \xi\to0, k\in \Z.
\ee
For $\ell=r,\ldots,m$, we take $u_\ell\in \lp{0}$ such that
\be \label{uell}
\wh{u_\ell}(\xi)=(i\xi)^\ell+\bo(|\xi|^{m+1}),\qquad \xi \to0,\; \ell=r,\ldots,m.
\ee
Now we define a vector function $h:=(h_1,h_2,\ldots,h_{r}):=(h_1,\theta_1,\ldots,\theta_{r-1})^\tp$, where $h_\ell:=\theta_{\ell-1}$ for $\ell=2,\ldots,r$ and
\be \label{h1}
h_1:=\theta_0+\sum_{\ell=r}^{m}
u_\ell*\theta_\ell=
\theta_0+\sum_{\ell=r}^m \sum_{k\in \Z} u_\ell(k) \theta_\ell(\cdot-k).
\ee
By the Hermite interpolation property of $\theta$ in \eqref{theta:hint}, we see that $h$ is a Hermite interpolant of order $r$ satisfying
\be \label{h:int}
H(k)=\td(k) I_r,\qquad \forall\; k\in \Z \quad \mbox{with}\quad H:=[h,h',\ldots,h^{(r-1)}].
\ee
Let $\vgu_H\in \lrs{0}{1}{r}$ be a vector sequence satisfying \eqref{vguH}.
By the definition of the function $h_1$ in \eqref{h1}, we have $\wh{h_1}(\xi)=\wh{\theta_0}(\xi)+\sum_{\ell=r}^m \wh{u_\ell}(\xi)\wh{\theta_\ell}(\xi)$.
Note that all $\wh{u_\ell}$ are $2\pi$-periodic.  From \eqref{uell}, we deduce that
\[
\wh{h_1}(\xi+2\pi k)=\wh{\theta_0}(\xi+2\pi k)+
\sum_{\ell=r}^m \wh{u_\ell}(\xi)\wh{\theta_\ell}(\xi+2\pi k)=\wh{\theta_0}(\xi+2\pi k)+
\sum_{\ell=r}^m (i\xi)^\ell \wh{\theta_\ell}(\xi+2\pi k)+\bo(|\xi|^{m+1})
\]
as $\xi \to 0$. Consequently, since $\vgu_H$ satisfies \eqref{vguH}, we conclude from \eqref{theta:approx} that for all $k\in \Z$,
\begin{align*}
\wh{\vgu_H}(\xi) \wh{h}(\xi+2\pi k)
&=\wh{h_1}(\xi+2\pi k)+
\sum_{\ell=1}^{r-1} (i\xi)^{\ell} \wh{\theta_\ell}(\xi+2\pi k)\\
&=\wh{\theta_0}(\xi+2\pi k)+
\sum_{\ell=1}^m (i\xi)^\ell \wh{\theta_\ell}(\xi+2\pi k)+\bo(|\xi|^{m+1})\\
&=\sum_{\ell=0}^m (i\xi)^\ell \wh{\theta_\ell}(\xi+2\pi k)+\bo(|\xi|^{m+1})
=\td(k)+\bo(|\xi|^{m+1})
\end{align*}
as $\xi \to 0$. That is, we proved
\be \label{h:approx}
\wh{\vgu_H}(\xi)\wh{h}(\xi+2\pi k)=\td(k)+\bo(|\xi|^{m+1}),\qquad \xi \to 0, k\in \Z.
\ee
Note that the Hermite interpolant $h\in (\CH{m})^r$ of order $r$  is no longer refinable.
On the other hand, since $\wh{c_1}(0)=\cdots=\wh{c_r}(0)=1$, by \cite[Lemma~3.4]{han03} or \cite[Lemma~5.6.7]{hanbook}, there exist $d_1,\ldots,d_r\in \lp{0}$ such that
\be \label{d:seq}
|\wh{d_{\ell}}(\xi)|\ge 1/2 \quad \forall\; \xi\in \R \quad \mbox{and}\quad
\wh{d_{\ell}}(\xi)=1/\wh{c_{\ell}}(\xi)
+\bo(|\xi|^{m+1}),\quad \xi \to 0, \forall\, \ell=1,\ldots,r.
\ee
We define a compactly supported initial vector function $f\in (\CH{m})^r$ by
\be \label{hsd:initial:f}
f:=(d_1*f_1,\ldots,d_r*f_r)^\tp,\quad \mbox{that is},\quad
\wh{f}(\xi)=(\wh{d_1}(\xi)\wh{h_1}(\xi),\ldots,
\wh{d_r}(\xi)\wh{h_r}(\xi))^\tp,
\ee
where $d_\ell*f_\ell:=\sum_{k\in \Z} d_\ell(k) f_\ell(\cdot-k)$ for $\ell=1,\ldots,r$.
Now it is trivial to deduce from \eqref{h:approx} and \eqref{d:seq}
that \eqref{initialf} holds. That is, $f$ is an admissible initial vector function.
Define $f_n:=\cd_a^n f$ for $n\in \N$, where the refinement operator $\cd_a$ is defined in \eqref{cd}.

(2)$\imply$(1).
By item (2) and $f_n:=\cd_a^n f$, we have $\lim_{n\to \infty} \|f_n-\phi\|_{(\CH{m})^r}=0$ and by $m\ge r-1\ge 0$ we particularly have $\lim_{n\to \infty} \|F_n-\Phi\|_{(\CH{})^{r\times r}}=0$, where
\be \label{FPhin}
F_n:=[f_n, f_n',\ldots,f_n^{(r-1)}]
\quad\mbox{and}\quad
\Phi:=[\phi, \phi',\ldots,\phi^{(r-1)}].
\ee
Note that
\[
\wh{f_n}(\xi)=\wh{\cd_a^n f}(\xi)=\wh{a_n}(2^{-n}\xi)\wh{f}(2^{-n}\xi)
=\wh{a_n}(2^{-n}\xi) \wh{B}(2^{-n}\xi) \wh{h}(2^{-n}\xi)=\wh{b_n}(2^{-n}\xi) \wh{h}(2^{-n}\xi),
\]
where $B:=\mbox{diag}(d_1,\ldots,d_r)$ and $b_n:=a_n*B$ satisfies
\[
\wh{b_n}(\xi):=\wh{a_n}(\xi)\wh{B}(\xi)=\wh{a}(2^{n-1}\xi)\cdots \wh{a}(2\xi)\wh{a}(\xi) \mbox{diag}(\wh{d_1}(\xi),\ldots, \wh{d_r}(\xi)).
\]
Then $f_n=2^n \sum_{k\in \Z} b_n(k) h(2^n\cdot-k)$ and
\[
F_n=[f_n, f_n',\ldots, f_n^{(r-1)}]= \sum_{k\in \Z} 2^n b_n(k)  H(2^n\cdot-k) \sD^{-n},
\]
where $H:=[h,h',\ldots,h^{(r-1)}]$ as in \eqref{h:int} and $\sD:=\mbox{diag}(1,2^{-1},\ldots,2^{1-r})$ as in \eqref{hsd:wn}.
Due to the Hermite interpolation property in \eqref{h:int}, since both $B$ and $\sD$ are diagonal matrices,
we have
\[
F_n(2^{-n}k)=2^n b_n(k) \sD^{-n}
=2^n (a_n*B)(k)\sD^{-n}=
2^n ((a_n\sD^{-n})*B)(k)
\]
for all $k\in \Z$. Now it follows from
\[
\left\| 2^n ((a_n \sD^{-n})*B)(\cdot) -\Phi(2^{-n}\cdot) \right\|_{\lrs{\infty}{r}{r}}
=\| F_n-\Phi\|_{\lrs{\infty}{r}{r}}
\le \| F_n-\Phi\|_{(\CH{})^{r\times r}}
\]
and $\lim_{n\to \infty} \|F_n-\Phi\|_{(\CH{})^{r\times r}}=0$ that
\be \label{bn:eq1}
\lim_{n\to \infty} \left\|
2^n ((a_n \sD^{-n})*B)(\cdot) -\Phi(2^{-n}\cdot) \right\|_{\lrs{\infty}{r}{r}}=0.
\ee
Since $|\wh{d_\ell}(\xi)|\ge 1/2$ for all $\xi \in \R$ and $\ell=1,\ldots,r$, we can define sequences $u_\ell$ on $\Z$ by
\[
\wh{u_\ell}(\xi):=1/\wh{d_\ell}(\xi), \qquad \ell=1,\ldots,d.
\]
Moreover, $\wh{u_\ell}(0)=1/\wh{d_\ell}(0)=\wh{c_\ell}(0)=1$ and all the sequences $u_\ell$ must have exponential decay. Define $U:=\mbox{diag}(u_1,\ldots,u_r)$. By $B*U=\td I_r$, we have
\[
\left\|
2^n a_n \sD^{-n} -\Phi(2^{-n}\cdot)*U \right\|_{\lrs{\infty}{r}{r}}
\le \|U\|_{\lrs{1}{r}{r}}
\left\| 2^n ((a_n \sD^{-n})*B)(\cdot) -\Phi(2^{-n}\cdot) \right\|_{\lrs{\infty}{r}{r}}
\]
Using \eqref{bn:eq1} and the above inequality, we conclude that $\lim_{n\to \infty} \left\|
2^n a_n(\cdot) \sD^{-n} -\Phi(2^{-n}\cdot)*U \right\|_{\lrs{\infty}{r}{r}}=0$.
Note that all entries in $\Phi$ are compactly supported continuous functions and all sequences $u_\ell$ have exponential decay with $\wh{u_1}(0)=\cdots\wh{u_r}(0)=1$. Consequently, by \cref{lem:to0}, we must have
$\lim_{n\to \infty} \left\|
\Phi(2^{-n}\cdot)*U -\Phi(2^{-n}\cdot) \right\|_{\lrs{\infty}{r}{r}}=0$. Now by the triangle inequality
\[
\left\|2^n a_n \sD^{-n} -\Phi(2^{-n}\cdot)\right\|_{\lrs{\infty}{r}{r}}
\le \left\|
2^n a_n \sD^{-n} -\Phi(2^{-n}\cdot)*U \right\|_{\lrs{\infty}{r}{r}}+\left\|
\Phi(2^{-n}\cdot)*U -\Phi(2^{-n}\cdot) \right\|_{\lrs{\infty}{r}{r}},
\]
we conclude that
$\lim_{n\to \infty}
\left\|2^n a_n(\cdot) \sD^{-n} -\Phi(2^{-n}\cdot)\right\|_{\lrs{\infty}{r}{r}}=0$.
Since $\sd_a^n (\td I_r)=2^n a_n$,
this proves \eqref{sdn:phi} and the Hermite subdivision scheme associated with mask $a$ is convergent with limiting functions in $(\CH{m})^r$. Hence we proved (2)$\imply$(1) without using the stability condition in \eqref{stability}.

(1)$\imply$(2). The condition in \eqref{stability} obviously implies \eqref{cond:phi}. Hence, by \cref{thm:hsd:mask}, since item (1) holds, its basis vector function $\varphi\in (\CH{m})^r$ of the Hermite subdivision scheme associated with mask $a$ must be refinable with $e_\tp \wh{\varphi}(0)=1$ and items (1) and (2) of \cref{thm:hsd:mask} hold.
On the other hand, the refinable vector distribution $\phi$ satisfies the refinement equation and $\wh{\vgu_a}(0)\wh{\phi}(0)=1$. However, due to item (1) of \cref{thm:hsd:mask}, the compactly supported distributional solutions to the refinement equation $\phi=2\sum_{k\in \Z} a(k) \phi(2\cdot-k)$ must be unique, up to a multiplicative constant. Because $\wh{\vgu_a}(0)=e_1$, we have $e_1^\tp \wh{\varphi}(0)=1=e_1^\tp \wh{\phi}(0)$. Therefore, we must have $\phi=\varphi$.
Hence, by $\varphi\in (\CH{m})^r$,
we must have $\phi\in (\CH{m})^r$.
Thanking to the stability condition in \eqref{stability} and $\phi\in (\CH{m})^r$, we conclude from \cite[Theorem~4.3]{han03} or
\cite[Corollary~5.6.12]{hanbook} that item (2) holds. This proves (1)$\imply$(2). This proves all the claims.
\ep

\end{document}